\documentclass[11pt, reqno]{amsart}

\usepackage[dvipsnames,usenames]{color}
\usepackage[colorlinks=true, urlcolor=NavyBlue, linkcolor=NavyBlue, citecolor=NavyBlue]{hyperref}
\usepackage{float}
\usepackage{amsmath, amsthm, amssymb}
\usepackage{amsfonts}
\usepackage{enumerate}
\usepackage{epsf}
\usepackage{psfrag}
\usepackage{amsmath}
\usepackage[all]{xy}
\usepackage[dvips]{graphicx}
\usepackage{epsfig}
\usepackage{epstopdf}
\usepackage{ifpdf}
\usepackage[UKenglish]{babel}
\usepackage{subfigure}
\usepackage{hyperref}
\usepackage{color}
\usepackage{marginnote}

\hypersetup{
linkbordercolor={1 0 0}, 
citebordercolor={0 1 0} 
}
\newtheorem{thm}{Theorem}[section]
\newtheorem{lemma}[thm]{Lemma}

\newtheorem{cor}[thm]{Corollary}

\newtheorem{conj}[thm]{Conjecture}

 \theoremstyle{definition}

\theoremstyle{remark}
\newtheorem{rmk}[thm]{Remark}

\newcommand{\rank}{\textup{rk}}

\numberwithin{equation}{subsection}
\numberwithin{enumi}{subsection}

\makeatletter

\newcommand{\Rmnum}[1]{\expandafter\@slowromancap\romannumeral #1@}
\makeatother

\author[Faramarz Vafaee]{Faramarz Vafaee}
\thanks{}
\address {Department of Mathematics, Michigan State University, East Lansing, MI 48824}
\email{vafaeefa@msu.edu}

\begin{document}
\title{On the knot Floer homology of twisted torus knots}
\maketitle
\begin{abstract}
In this paper we study the knot Floer homology of a subfamily of twisted $(p, q)$ torus knots where $q \equiv\pm1 \pmod{p}$. Specifically, we classify the knots in this subfamily that admit L-space surgeries. To do calculations, we use the fact that these knots are $(1, 1)$ knots and, therefore, admit a genus one Heegaard diagram.
\end{abstract}
\section{Introduction}\label{section:1}
Heegaard Floer theory consists of a set of invariants of three- and four-dimensional manifolds \cite{Ozsvath2004a}. For $Y$ a closed three manifold, one example of such invariants is $\widehat{HF}(Y)$, which is a Spin$^c$-graded abelian group where the Spin$^c$ structures of $Y$ are in one to one correspondence with the elements of $H^2(Y; \mathbb{Z})$. Lens spaces have the simplest Heegaard Floer homology, that is, $\widehat{HF}(Y, \mathfrak{s}) \cong \mathbb{Z}$ for each $\mathfrak{s}$ in Spin$^c$(Y). By definition, a rational homology three-sphere with the same property is called an \emph{L-space}.

A knot $K \subset S^3$ is called an \emph{L-space knot} if performing $n$-surgery on $K$ results in an L-space for some positive integer $n$. Any knot with a positive lens space surgery is then an L-space knot. In \cite{Berge}, Berge gave a conjecturally complete list of knots that admit lens space surgeries including torus knots \cite{Moser1971}. Therefore it is natural to look beyond Berge's list for L-space knots. Examples include the $(-2, 3, 2n+1)$ pretzel knots (for positive integers n) \cite{Paper1996, FS1980, Ath}, which are known to live outside of Berge's collection when $n \ge 5$ \cite{Mattman2000}. It is also proved in \cite{Lidman} that these 3-strand pretzel knots are the only pretzel knots with L-space surgeries. Another source of L-space knots is within the set of cable knots. By combining work of Hedden \cite{Hedden2009} and Hom \cite{Hom2011a}, the $(p, q)$ cable of a knot, $K$, is an L-space knot if and only if $K$ is an L-space knot and $q/p \ge 2g(K) - 1$. 

\begin{figure}[t]
\begin{center}
\psfrag{t}{$\tau$}   
\psfrag{p}{$(p, q)$}
\psfrag{r}{\small{torus knot}}
\psfrag{q}{\tiny{$r$}}
\psfrag{u}{\tiny{full twists}}
\psfrag{s}{\tiny{on}}
\psfrag{v}{\tiny{$s$ strands}}
 \includegraphics[scale=.45]{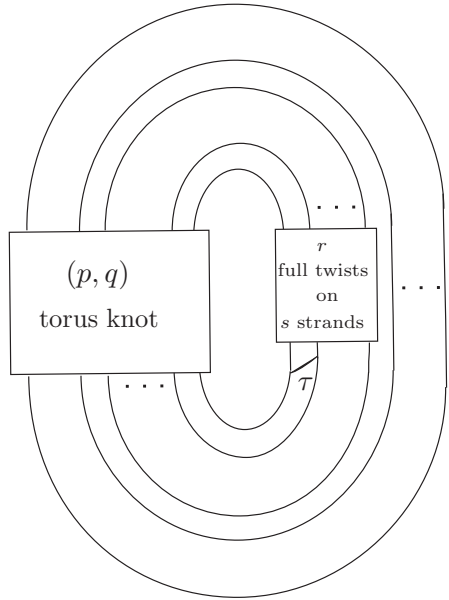}
       \caption{\small{A $(p,q)$ torus knot with $r$ positive full twists on $s$ adjacent strands. (Here, $p$ denotes the longitudinal winding.) The arc $\tau$ is a one-bridge, i.e. it divides the knot into two arcs, where one arc is unknotted and the other arc can be trivialized (unknotted) by sliding one or both of its endpoints along the \emph{a priori} unknotted arc. In order to make sense of adjacency of strands, we need to have the standard presentation of a torus knot. Note that where the twist occurs is irrelevant. }} 
\label{fig1}
\end{center}
\end{figure}

The primary purpose of this paper is to investigate L-space knots in the family of \emph{twisted torus knots}, $K(p, q;s, r)$, which are defined to be $(p, q)$ torus knots with $r$ full twists on $s$ adjacent strands where $\displaystyle 0< s < p$. See Figure \ref{fig1}. Watson proved in \cite{Watson} that the knots $K(3, 3k+2;2, 1)$ are L-space knots ($k > 0$). We generalize this result in Corollary~\ref{cor3} by showing that  all twisted $(3, q)$ torus knots admit L-space surgeries ($q$ is a positive integer that does not divide $3$). 

To the best of our knowledge, the examples mentioned are the only previously known explicit families of L-space knots. If $K$ is a quasi-alternating knot with unknotting number one, then the preimage of an unknotting arc in the branched double cover of $K$ is a knot in an L-space with an $S^3$ surgery (see \cite{montesinos1973}, \cite[Section 8.3]{ozsvath2005knots}, and \cite[Proposition 3.3]{Ozsvath2005}). The dual to this curve is therefore a knot in $S^3$ with an L-space surgery, so either it or its mirror image is an L-space knot. However, at present, there is no explicit parametrization of the knots that arise in this way. In this paper, we classify all the L-space twisted $(p, q)$ torus knots with $q = kp \pm 1$. The question of what happens when $q \ne kp \pm 1$ remains unanswered. Our examples include the L-space pretzel knots as a proper subfamily since the $(-2, 3, 2m+3)$ pretzel knot is isotopic to $K(3, 4; 2, m)$ for $m\ge1$.

We now state the main result of the paper. With the above notation:
\\

{\newtheorem*{theorem*}{Theorem} \label{theorem:1}\noindent \textbf{Theorem 1.} {\it For $p \ge 2$, $k \ge 1$, $r > 0$ and $0 < s <p$, the twisted torus knot, $K(p, k p \pm 1; s, r)$, is an L-space knot if and only if either $s = p - 1$ or $s \in \left \{2, p - 2 \right \}$ and $r = 1$.}} 
\\

A key ingredient of the proof is the observation that all of the twisted torus knots being studied are \emph{$(1, 1)$ knots}, that is, knots that can be placed in one-bridge position with respect to a genus one Heegaard splitting of $S^3$. Thus, the knot is comprised of two properly embedded unknotted arcs, one in each of the two solid tori of the Heegaard splitting. These arcs meet along their endpoints so that their union is equal to the knot. 

From the perspective of knot Floer homology, (1,1) knots are particularly appealing. It was first observed by Goda, Morifuji, and Matsuda \cite{Goda2005} that $(1,1)$ knots are exactly those knots that can be presented by a \emph{doubly-pointed Heegaard diagram} of genus one. The chain complex for knot Floer homology is defined in terms of a doubly-pointed Heegaard diagram. As shown by Ozsv\'ath and Szab\'o \cite{Ozsvath2004}, for knots admitting a genus one diagram, knot Floer homology can be computed combinatorially and efficiently. 

The outline of the paper is as follows: Section \ref{section:2} introduces the theory of $(1, 1)$ knots and presents how to draw a genus one Heegaard diagram for $(1, 1)$ knots via an explicit example. Section \ref{section:3} contains the main result of the paper, as well as the corollaries. In the final section, we state some questions that address future research.
\\

\noindent {\bf Acknowledgements.} I would like to express my sincere gratitude to Matthew Hedden for suggesting this project to me and for his invaluable guidance as an advisor. I would also like to thank Adam Giambrone for his detailed and thoughtful comments on an earlier draft of this paper, David Krcatovich for numerous enlightening and instructive discussions, and Allison Moore for some helpful email correspondence and her interest in my work. Finally, I am grateful to the anonymous referee for advantageous suggestions and favorable comments.  

\section{Background and preliminary lemmas}\label{section:2}

We start this section by showing that the knots $K(p, kp \pm 1; s, r)$ are $(1, 1)$ knots. Next, we explain an algorithm which produces genus one Heegaard diagrams for the twisted torus knots with a $(1, 1)$ decomposition. Finally, we assemble some preliminary facts needed to prove Theorem~\ref{theorem:1}.

\begin{figure}[t!]
 \centering
 \subfigure[\label{fig2:subfig1}]{
\psfrag{t}{$\tau$}
  \includegraphics[scale=.45]{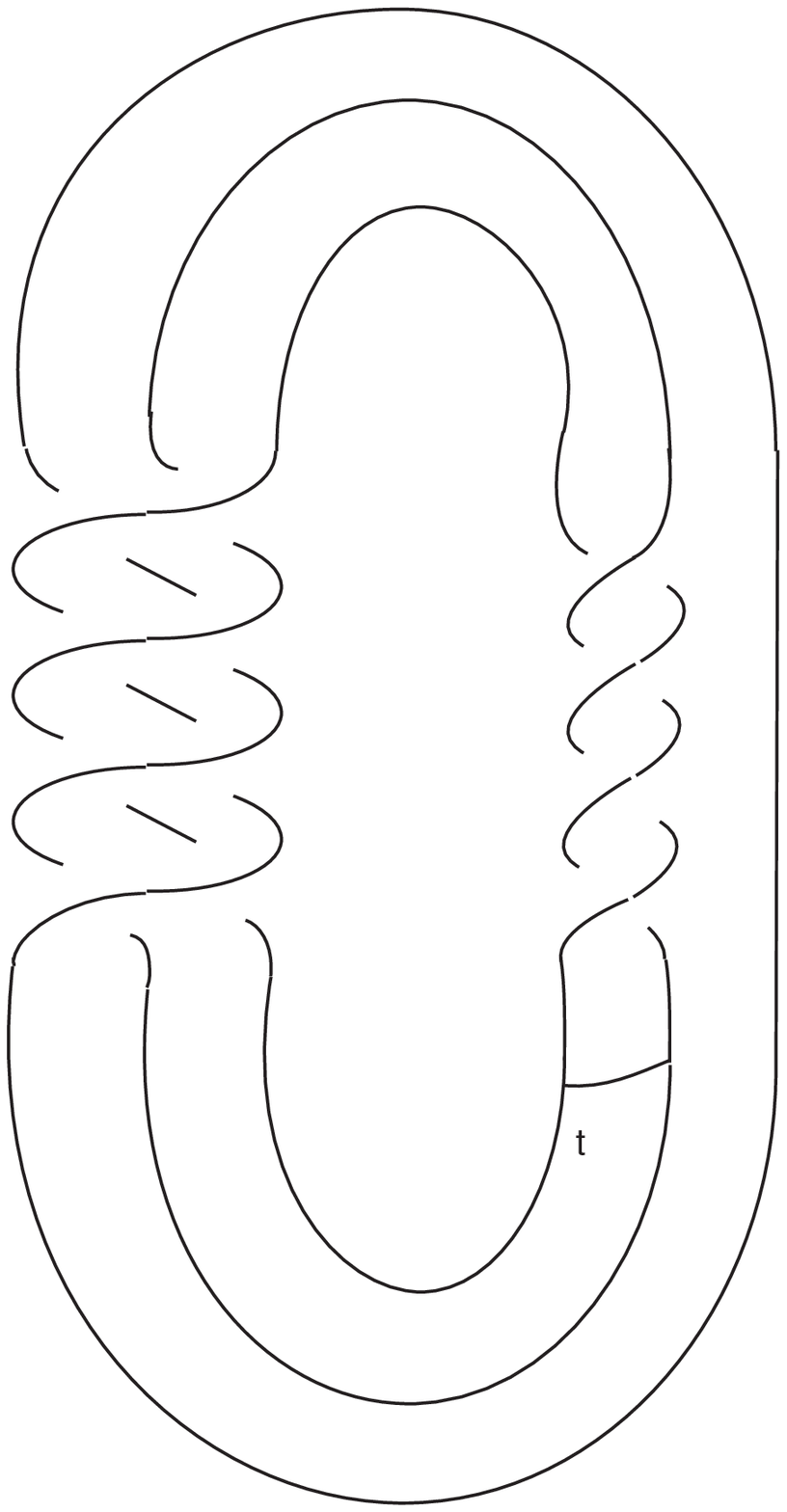}
  
   }
 \subfigure[]{
\psfrag{t}{$\tau$}
  \includegraphics[scale=.45]{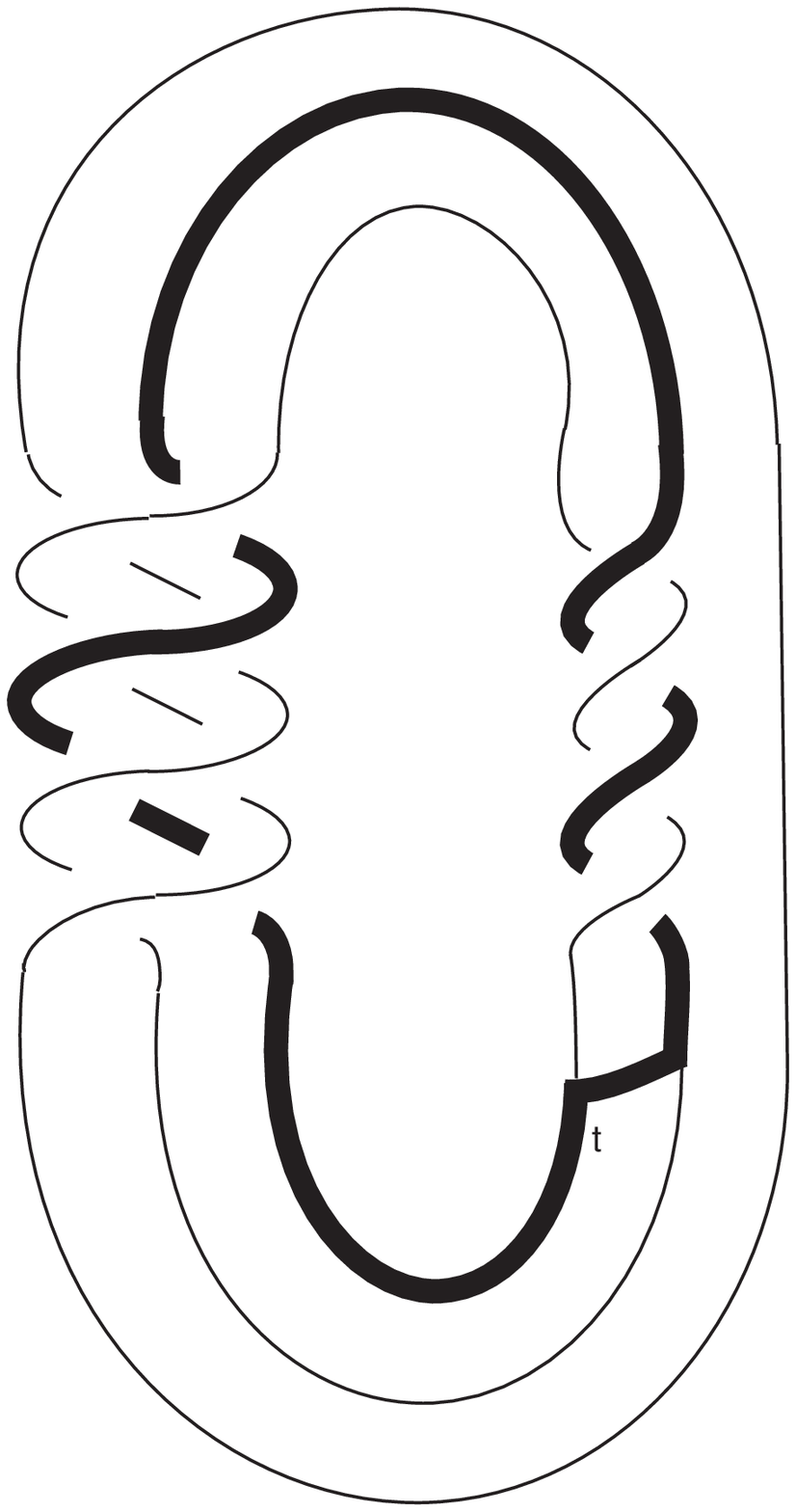}
  \label{fig2:subfig2}
   }
       \caption{\small{A $(3, 4)$ torus knot with two positive full twists on two adjacent strands. The one-bridge is indicated by $\tau$.}}
\label{fig2}
\end{figure}

\subsection{(1, 1) knots and genus one Heegaard diagrams}\label{subsec1}Let $p$ and $q$ be relatively prime positive integers and let $r$ and $s$ be integers. We denote the knot illustrated in Figure \ref{fig1} by $K(p, q;s, r)$.  Let $\tau$ be the arc indicated in Figure~\ref{fig1}. By untying the crossings of the $r$ full twists above the arc through edge slides along the arc, we will show that $\tau$ becomes a one-bridge for $K(p, q;s, r)$ provided that $q=kp\pm1$. See Figure \ref{fig2} for an explicit example. It has been a long standing question of whether or not any twisted torus knot, with $q$ that is not of the form $kp \pm 1$, is a $(1, 1)$ knot. In 1991, Morimoto, Sakuma, and Yokota conjectured that the answer is negative:

{\conj[\cite{Morimoto2008}, Conjecture 1.3] $K(p, q;2, r)$ admits no $(1, 1)$ decomposition unless either $p \equiv \pm 1 \pmod {q}$, or $q \equiv \pm 1 \pmod{p}$, or $r = 0, \pm1$.}
\\

 Having $s = 2$ does not seem to play an important role in the conjecture and, in fact, we expect a similar conjecture to hold when the twisting is on any number of strands. Bowman, Taylor, and Zupan have proved this conjecture when the number of twists is large \cite[Theorem 1.1]{bowman2014}.

In the rest of this subsection, we give an explicit construction of a genus one doubly-pointed Heegaard diagram via a specific example, namely $K= K(3, 4;2, 2)$. See Figure \ref{fig2}. This example should help clarify the strategy we use for our calculations.  

\begin{figure}[H]
 \begin{center}
 \subfigure[   \label{fig3:subfig1}]{
  \psfrag{z}{$z$}  
\psfrag{w}{$w$}
  \includegraphics[scale=.4]{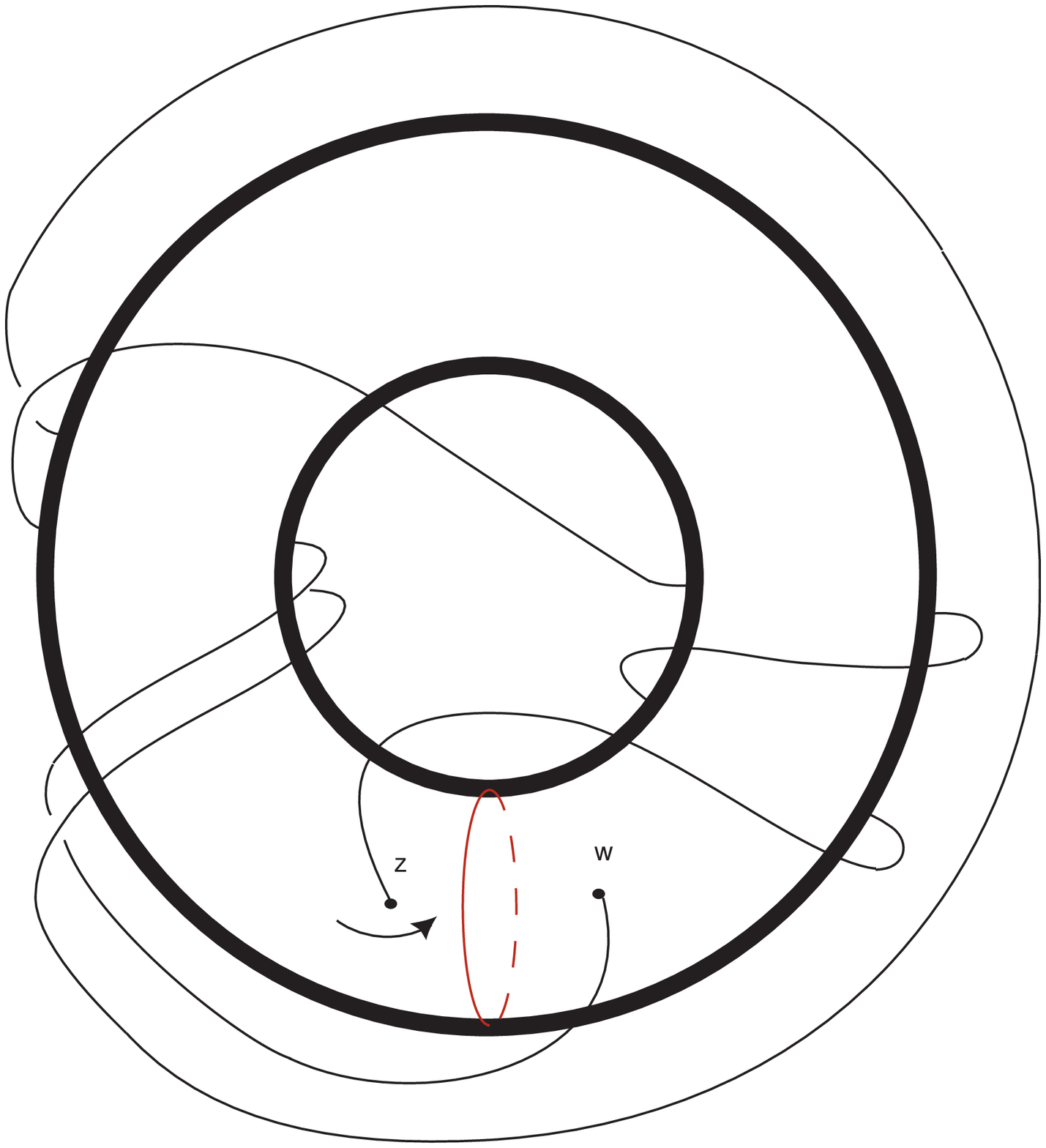}

   } 
\subfigure[   \label{fig3:subfig2}]{
  \psfrag{z}{$z$}  
\psfrag{w}{$w$}
  \includegraphics[scale=.4]{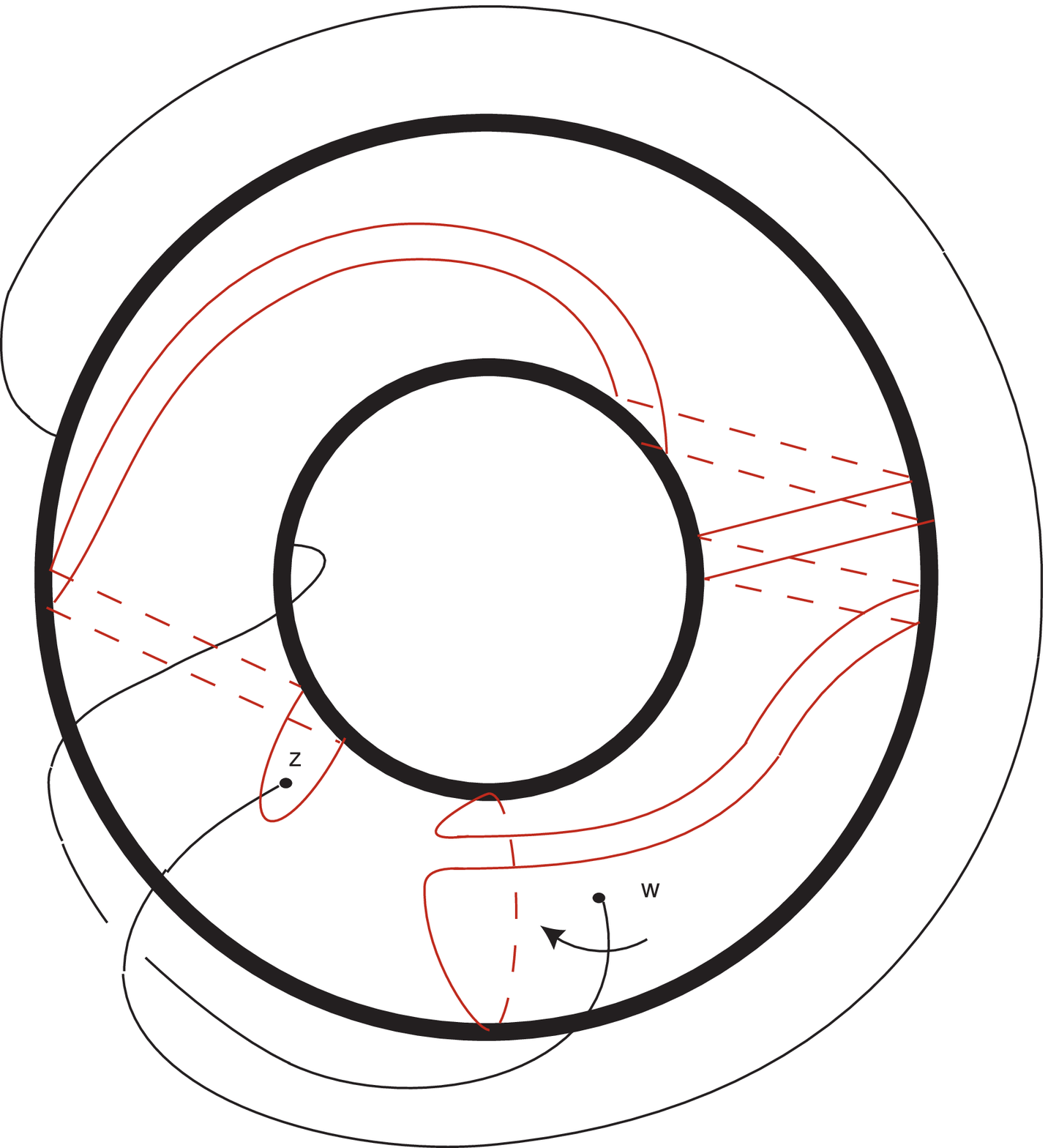}

   }
 \subfigure[   \label{fig3:subfig3}]{
  \psfrag{z}{$z$}  
\psfrag{w}{$w$}
  \includegraphics[scale=.38]{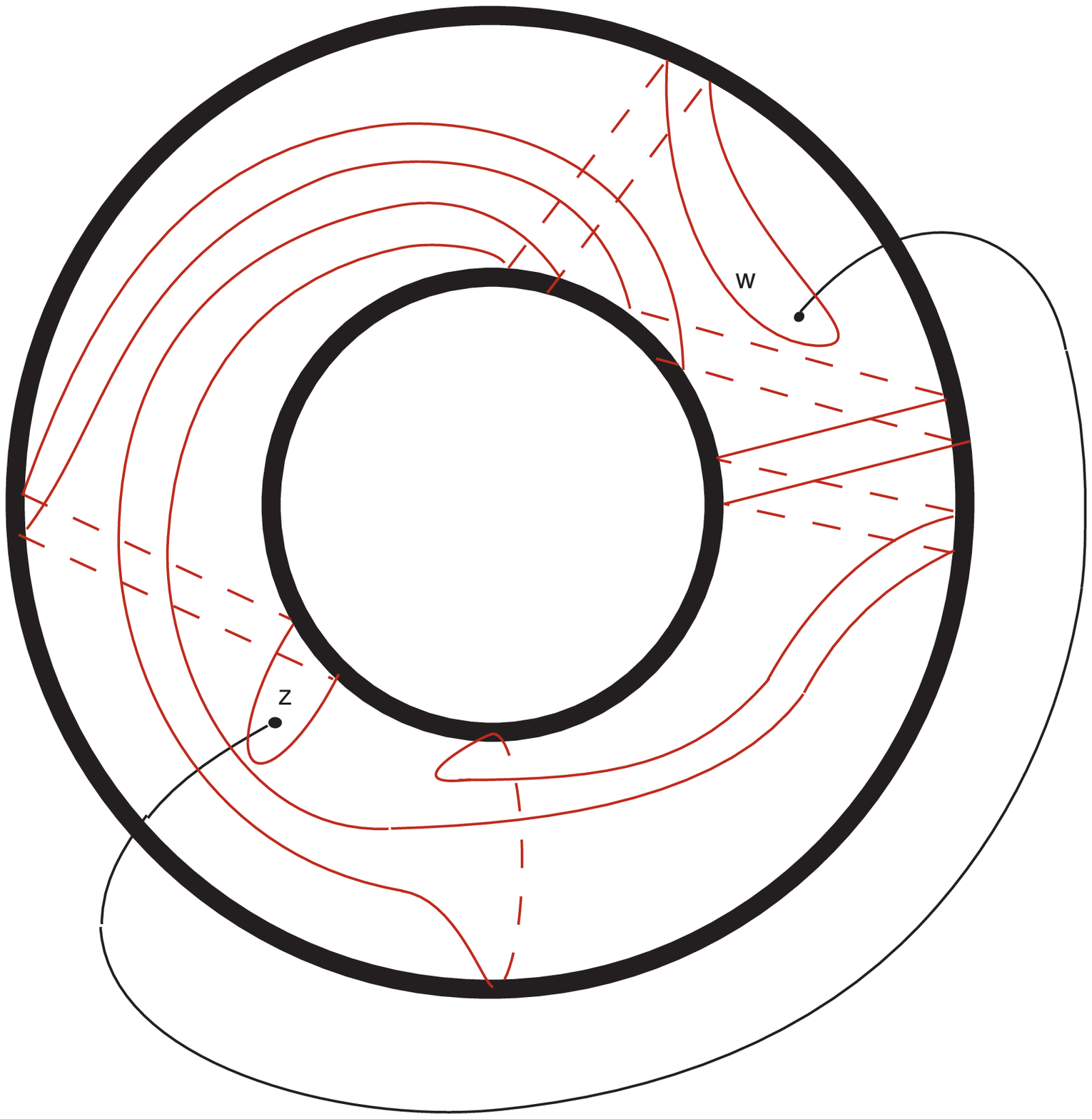}

   }
 \subfigure[   \label{fig3:subfig4}]{
  \psfrag{z}{$w$}  
\psfrag{w}{$z$}
\psfrag{1}{$1$}
\psfrag{2}{$2$}
\psfrag{3}{$3$}
\psfrag{4}{$4$}
\psfrag{5}{$5$}
\psfrag{6}{$6$}
\psfrag{7}{$7$}
\psfrag{8}{$8$}
\psfrag{9}{$9$} 
  \includegraphics[scale=.34]{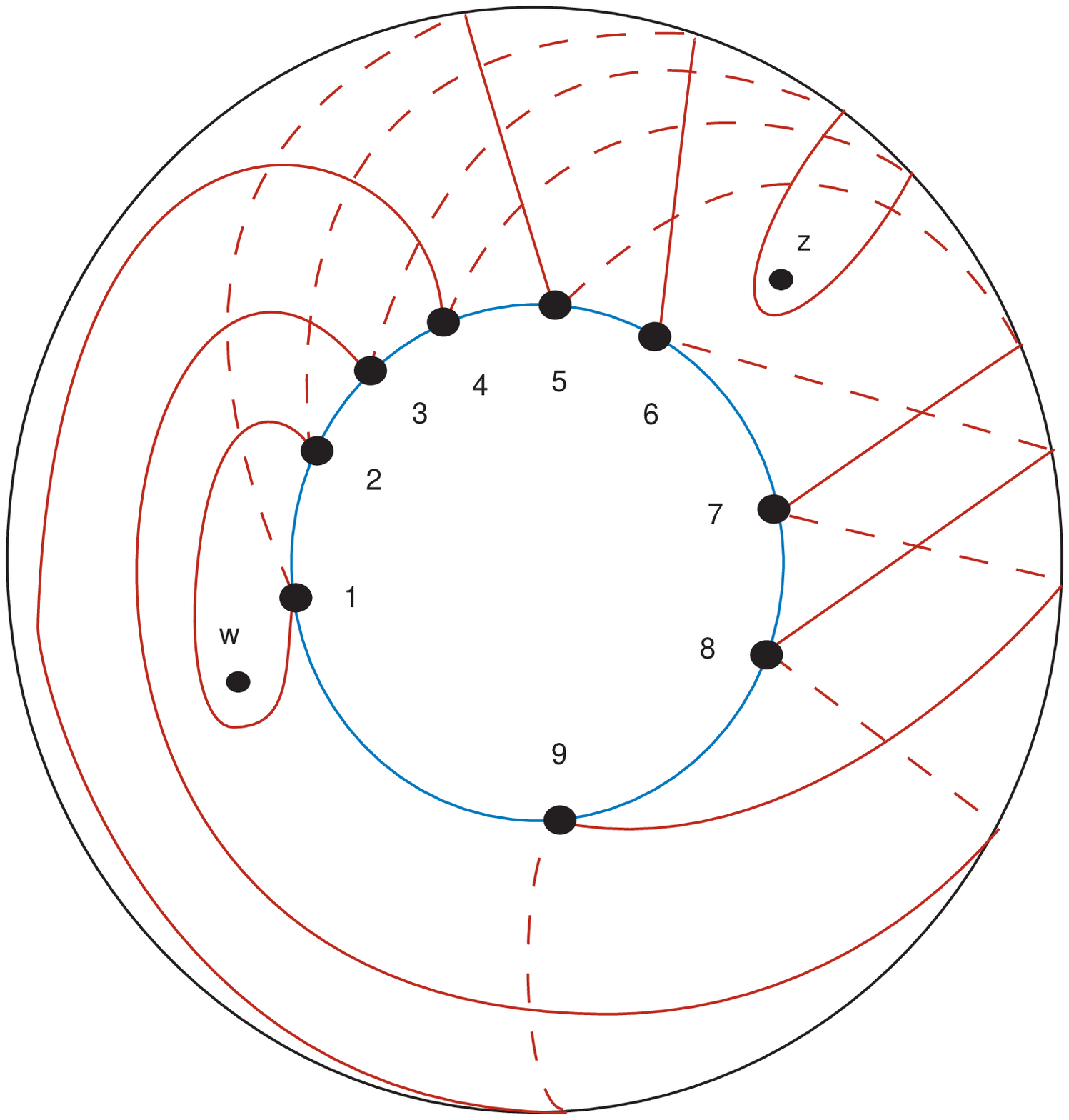}

   }
\caption{\small{The process of obtaining a genus one Heegaard diagram for the $(3, 4)$ torus knot with two positive full twists on two adjacent strands. In the algorithm discussed in Section~\ref{subsec1}, Figure~\ref{fig3:subfig2} corresponds to Step~1, and also Figure~\ref{fig3:subfig3} corresponds to, simultaneously, implementing Step~2 and Step~3. Note that the torus (in bold) corresponds to a neighborhood of the bold curve of Figure \ref{fig2:subfig2}. Note also that the $\alpha$ curve is drawn in red and the $\beta$ curve is drawn in blue.}}{%
 }
\label{fig3}{%
}
\end{center}
\end{figure}

We now describe a procedure to see that the arc $\tau$ (indicated in Figure \ref{fig2}) is a \emph{one-bridge}, i.e. it divides the knot $K$ into two arcs, where one arc is \emph{a priori} unknotted and the other arc can be trivialized (unknotted) by sliding one or both endpoints of this arc along the bold curve in Figure~\ref{fig2:subfig2}. (See \cite{ording2006} for a detailed discussion on how to produce a genus one Heegaard diagram for a certain family of $(1, 1)$ knots.) The closed curve indicated in bold is the union of the one-bridge, $\tau$, and the \emph{a priori} unknotted arc. Therefore, its neighborhood is an unknotted torus. In Figure~\ref{fig3} we show, diagrammatically, how to use the one-bridge and the unknotting process to obtain a Heegaard diagram for the knot $K$. (The red and blue curves in Figure~\ref{fig3} ($\alpha$ and $\beta$ curves respectively) are the boundaries of the meridional disks corresponding to the two solid tori of the genus one Heegaard splitting of $S^3$.) We do this by trivializing the arc living in the complement of the torus. To begin, move the $z$ base point in the counterclockwise direction, making sure that the $z$ base point passes to the left of the $w$ base point, as otherwise we would create more crossings rather than simplify the arc. See Figure~\ref{fig3:subfig2}.  Now move the $w$ base point in the clockwise direction, passing to the left of the $z$ base point. See Figure~\ref{fig3:subfig3}. That completes the construction of the genus one Heegaard diagram. See Figure~\ref{fig3:subfig4}. 

This construction can be generalized to an algorithm with three steps to produce a genus one Heegaard diagram for $K(p, kp\pm1; s, r)$. Note that the number of longitudinal and meridional windings is dictated by the arc living in the torus complement:
\\

\noindent Step 1: Wind the $z$ base point once around the torus in the counter clockwise direction. Note that $z$ traverses the torus $(k + r)$ times meridionally.

\noindent Step 2: Wind the $w$ base point $(s -2)$ times in the clockwise direction. Note that each time $w$ traverses the torus $(k + r)$ times meridionally.

\noindent Step 3: Finally, wind the $w$ base point $(p - s)$ times, longitudinally, to completely trivialize the arc (in the sense that the planar projection of the arc no longer has any self-intersection). Note that each longitudinal winding goes through $k$ meridional moves. 

{\rmk \label{rmk} To trivialize the part of the knot that lives outside of the torus, we isotope the base points, $z$ and $w$, on the torus which forces the $\alpha$ curve to be perturbed. Specifically, in a neighborhood of the base points, the isotopy drags one (or more) sub-arc(s) of $\alpha$.}
\\

\noindent Note that the Heegaard diagram in Figure~\ref{fig3:subfig4} may be represented by a rectangle with canonical identification implicit. See Figure~\ref{fig5:subfig1}.
\subsection{Lifted Heegaard diagrams, L-space knots, and {\bf$CFK^-$}}\label{subsec2} For $K \subset S^3$ a knot, let $CFK^-(K)$ denote the knot Floer complex associated to $K$ \cite{Ozsvath2004}. Fortunately, computing $CFK^{-}(K)$ for a $(1, 1)$~knot $K$ is purely combinatorial. We refer the interested reader to \cite[p.89]{Ozsvath2004} and \cite{Goda2005} for further details. To analyze holomorphic disks in the torus, it is convenient to pass to the universal covering space $\pi : \mathbb{C} \rightarrow T$. Given the base points $z$ and $w$ in $T$, $\pi^{-1}(z)$ and $\pi^{-1}(w)$ lift to affine lattices $Z$ and $W$, respectively. Also let $\{ \tilde{\alpha}_i \}$ and $ \{ \tilde{\beta}_j \}$ be the connected components of $\pi^{-1}(\alpha)$ and $\pi^{-1}(\beta)$, respectively. Now, given two intersection points $x$ and $y$ between $\alpha$ and $\beta$, the element $\phi \in \pi_{2}(x, y)$ is a Whitney disk that has Maslov index one and admits a holomorphic representative if and only if there is a bigon $\tilde{\phi} \in \pi_2(\tilde{x}, \tilde{y})$ with Maslov index one, where $\tilde{x}$ and $\tilde{y}$ are lifts of $x$ and $y$, intersection points between  $\tilde{\alpha}_i$ and $\tilde{\beta}_j$ (for some $i$ and $j$). In particular, $\mathcal{M}(\tilde{\phi}) \cong \mathcal{M}(\phi)$. See \cite{Ozsvath2004} for the notation 
\newpage
\newpage
\begin{figure}[H]
 \centering
 \subfigure[   \label{fig140:subfig1}]{
  \psfrag{z}{$z$}  
\psfrag{w}{$w$}
  \includegraphics[scale=.33]{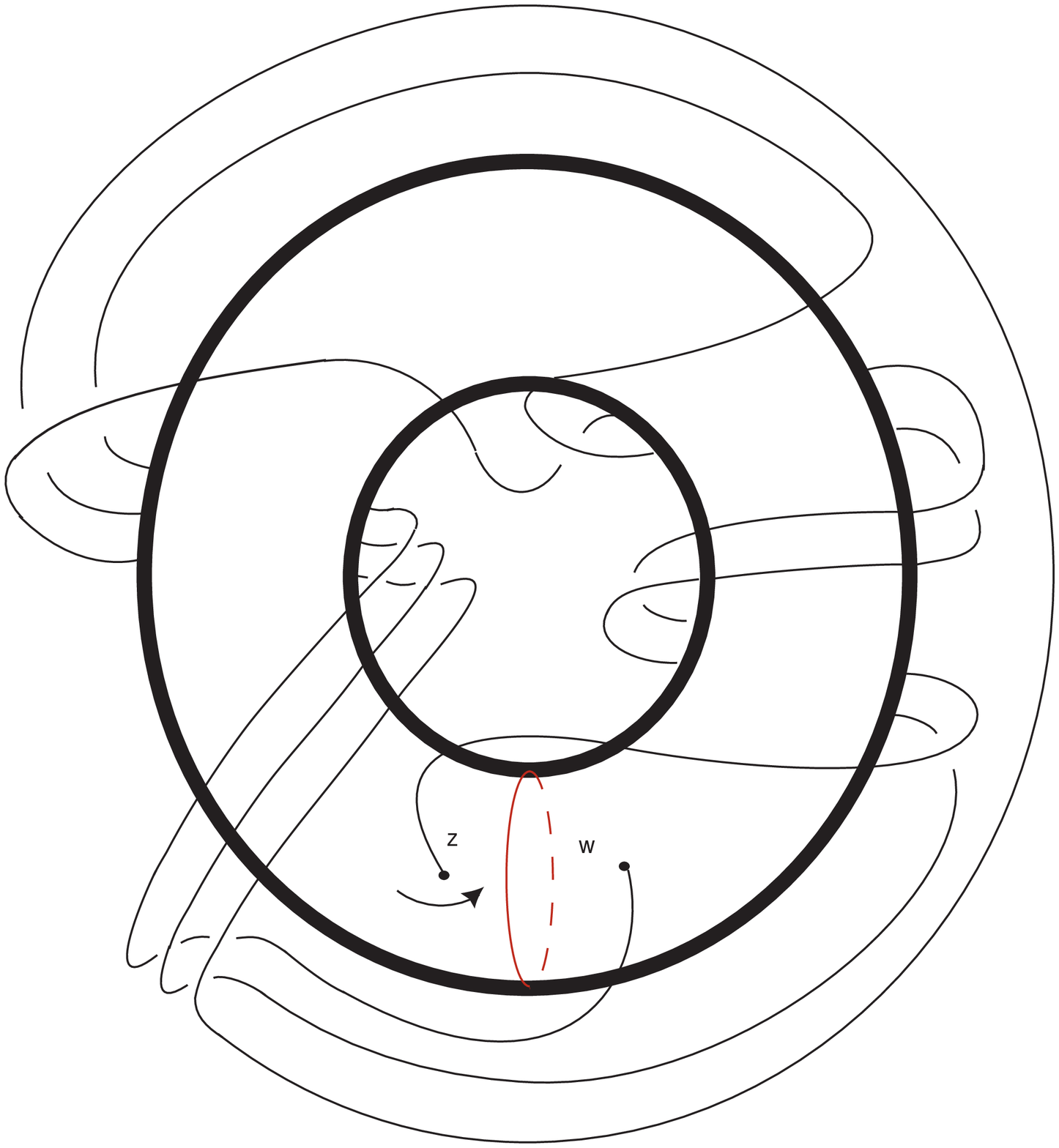}

   } 
\subfigure[   \label{fig140:subfig2}]{
  \psfrag{z}{$z$}  
\psfrag{w}{$w$}
  \includegraphics[scale=.33]{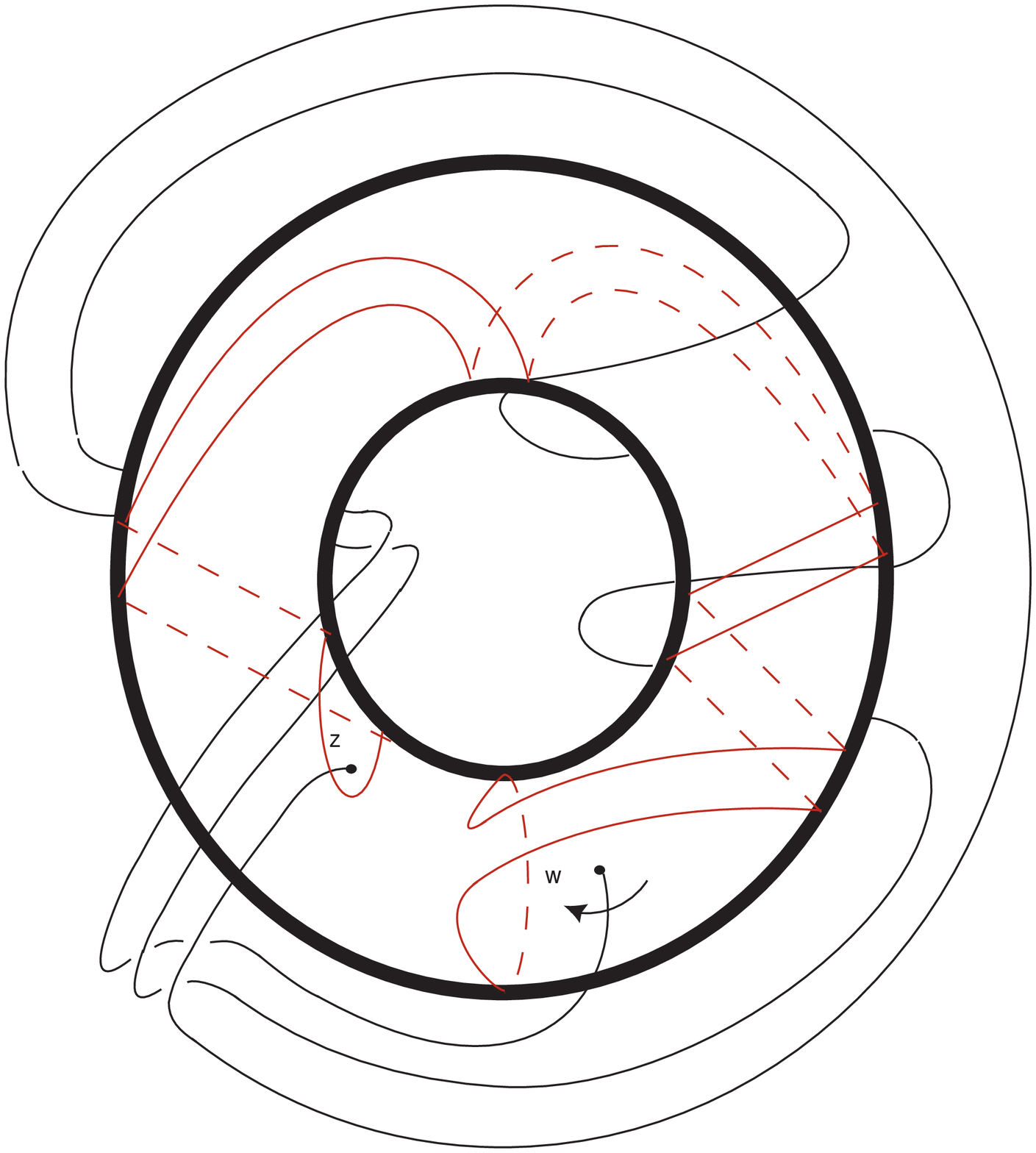}

   }
 \subfigure[   \label{fig140:subfig3}]{
  \psfrag{z}{$z$}  
\psfrag{w}{$w$}
  \includegraphics[scale=.37]{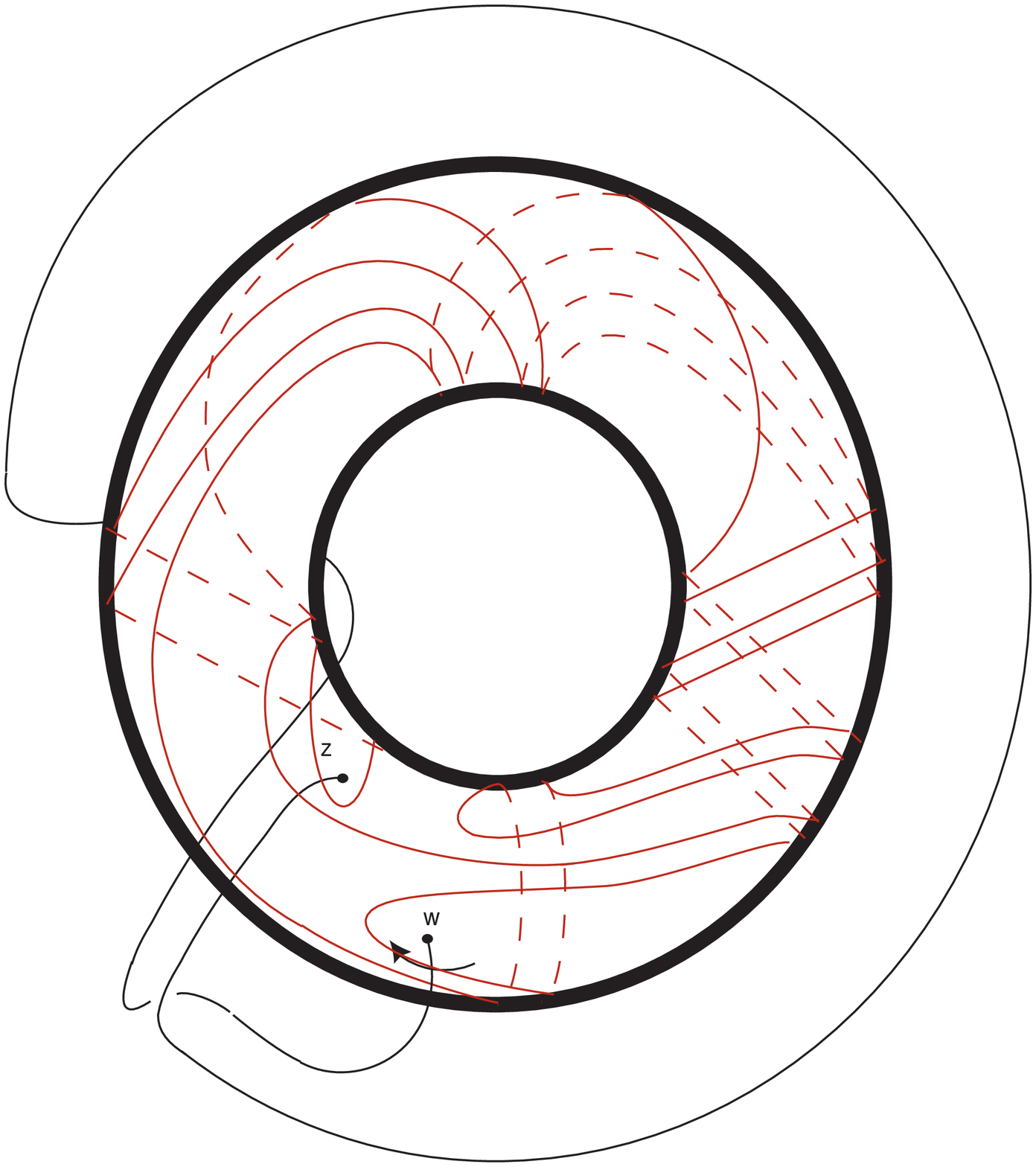}

   }
 \subfigure[   \label{fig140:subfig4}]{
  \psfrag{z}{$z$}  
\psfrag{w}{$w$}
  \includegraphics[scale=.37]{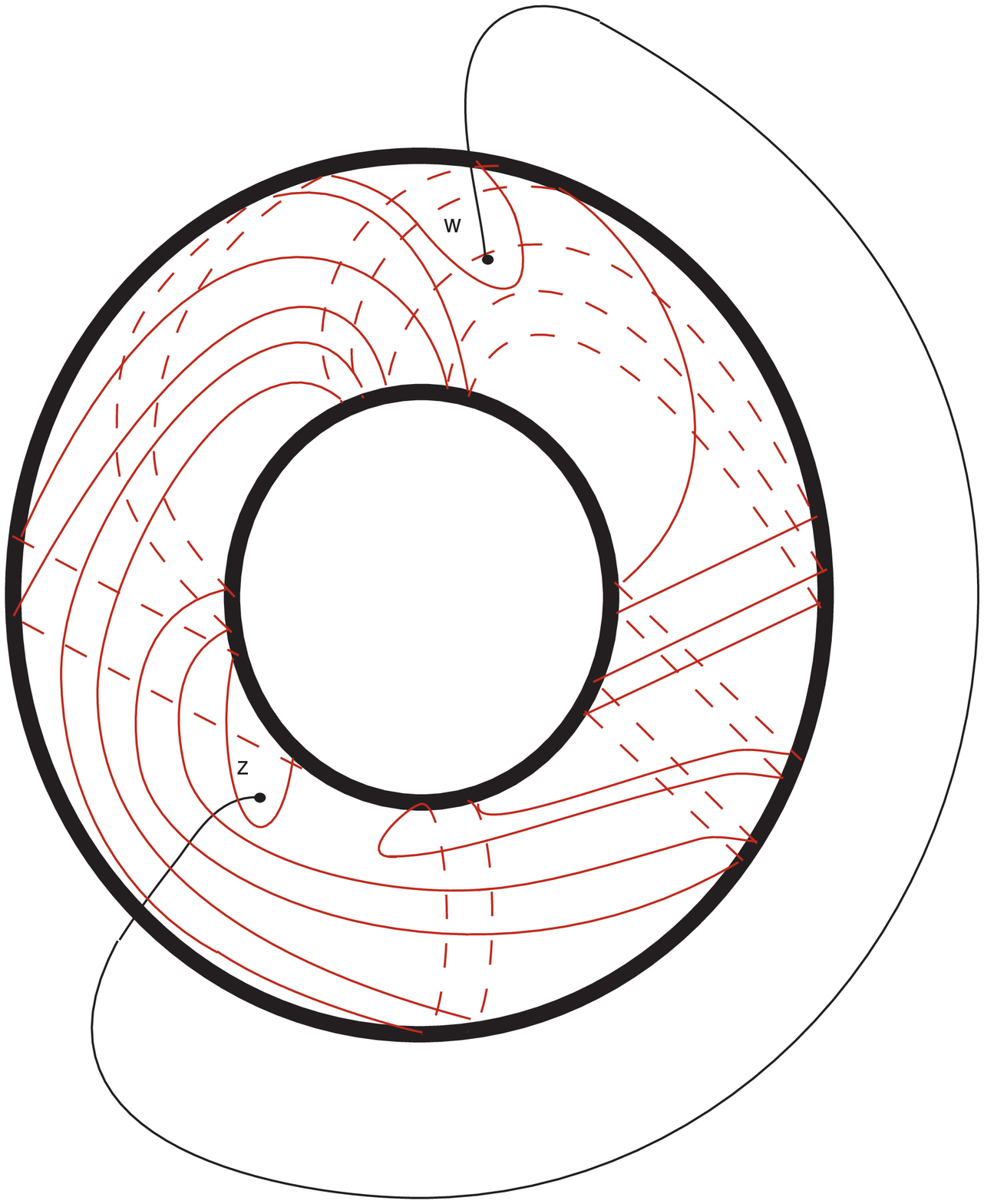}

   }

\caption{\small{The process of drawing a genus one Heegaard diagram for the $(4, 5)$ torus knot with two positive full twists on three adjacent strands. Figure~\ref{fig140:subfig2}, Figure~\ref{fig140:subfig3}, and Figure~\ref{fig140:subfig4} correspond to Step 1, Step 2, and Step 3, respectively, in the algorithm discussed in Section~\ref{subsec1}. The $\alpha$ curve is drawn in red. The base points must pass to the left of each other, as otherwise we would create more crossings rather than simplify the arc living in the torus complement.}}{%
 }
\label{fig14}{%
}
\end{figure}
\noindent 
\newpage 
\noindent 
used above. Figure~\ref{fig7:subfig2} shows a Heegaard diagram for $K = K(3, 4; 2, 2)$ that has been lifted to $\mathbb{C}$. Also, Figure \ref{fig8} represents $CFK^-(K)$. An  L-space knot $K$ can be thought of as a knot with the simplest knot Floer invariants. To make sense of this fact, note that \cite{Ozsvath2004a}
\begin{equation}\label{eqn:t}
\Delta_K(T) = \sum_{m, \mathfrak{s}} (-1)^m \text{rk } \widehat{HFK}_m(K, \mathfrak{s}) T^{\mathfrak{s}},
\end{equation}
where $\Delta_K(T)$ is the symmetrized Alexander polynomial of $K$. We observe that the total rank of $\widehat{HFK}(K)$ is bounded below by the sum of the absolute values of the coefficients of the Alexander polynomial of $K$. A necessary condition for $K$ to be an L-space knot is for this bound to be sharp. The following lemma turns out to be useful during the course of proving Part (c) of Theorem~\ref{theorem:2}. See \cite[Theorem 1.2]{Ath} for the complete statement. 

{\lemma \label{lem0}
 Assume that $K\subset S^3$ is a knot for which there is an integer $p$ such that $S^3_p(K)$ is an L-space. Then
\[
\begin{array}{cc}
\rank \text{ }\widehat{HFK}(K, \mathfrak{s}) \le 1 & \forall \mathfrak{s} \in \mathbb{Z}.
\end{array}
\]
\noindent In particular, all of the non-zero coefficients of $\Delta_K(T)$ are $\pm 1$.
}
\\

\noindent Therefore, if the absolute value of one of the coefficients of $\Delta_K(T)$ is greater than one, then $K$ is not an L-space knot. We end this subsection by noting that a knot Floer complex with a staircase-shape (as in Figure \ref{fig8}) represents an L-space knot. Such a complex  has a basis $\{ x_1, x_2, ..., x_m \}$ for $CFK^{\infty}(K)$ (defined in \cite{Ozsvath2004}) such that 

\begin{equation}\label{staircase}
\begin{array}{ccc} 
\partial x_i = &x_{i-1}+x_{i+1} & \text{for $i$ even} \\
\partial x_i = &0                      & \text{otherwise,}
\end{array}
\end{equation}

\noindent where the arrow from $x_i$ to $x_{i-1}$ is horizontal and the arrow from $x_i$ to $x_{i+1}$ is vertical. (We refer the reader to \cite[Section 6]{Hom} for the concept of a knot Floer complex with a staircase-shape.) The following corollary is a consequence of \cite[Remark 6.6]{Hom}.
{\cor \label{lem1}For a knot $K \subset S^3$, if $CFK^-(K)$ has a staircase-shape, then $K$ is an L-space knot.}
\section{Proof of the main theorem}\label{section:3}
This section is devoted to the proof of the main result of the paper. For the sake of the proof, it will convenient to restate Theorem \ref{theorem:1} in the following equivalent form:
{\thm \label{theorem:2}For $p \ge 2$, $k \ge 1$, $r > 0$ and $0 < s <p$, we have that $K(p, k p \pm 1; s, r)$:

(a) is an L-space knot if $s = p - 1$,

(b) is an L-space knot if $r = 1$ and $s \in \left \{ 2, p - 2 \right \}$, and

(c) does not admit any L-space surgeries otherwise.} 
\begin{figure}[t!]
 \centering
 \subfigure[\small{A Heegaard diagram for the $(3, 4)$ torus knot with two positive full twists on two adjacent strands \label{fig5:subfig1}}]{
\psfrag{z}{$z$}  
\psfrag{w}{$w$}  
\psfrag{1}{$1$}
\psfrag{2}{$2$}
\psfrag{3}{$3$}
\psfrag{4}{$4$}
\psfrag{5}{$5$}
\psfrag{6}{$6$}
\psfrag{7}{$7$}
\psfrag{8}{$8$}
\psfrag{9}{$9$} 
  \includegraphics[scale=.35]{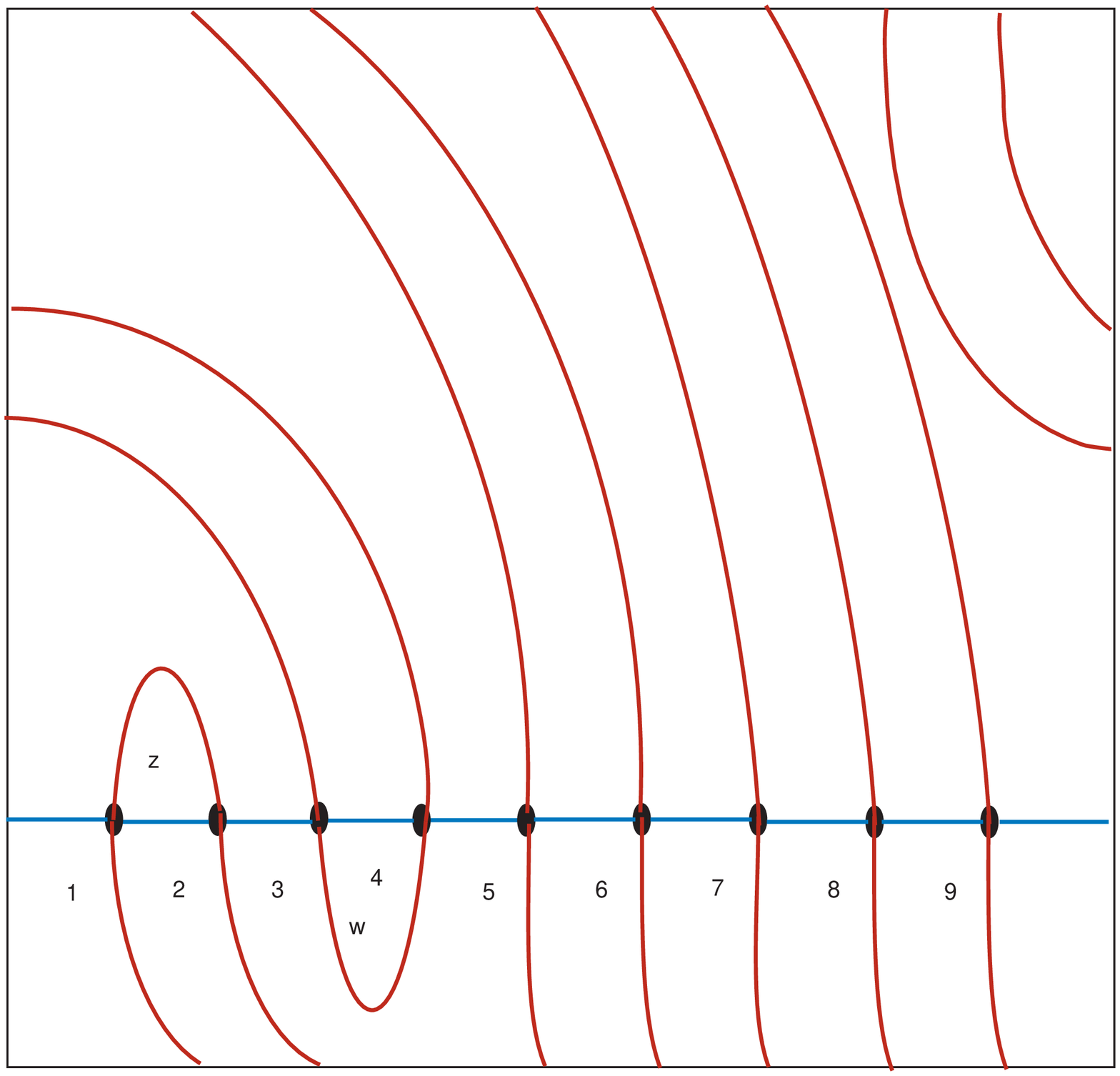}
   
   }
 \subfigure[\small{The general form of a Heegaard diagram for $K(p, kp \pm 1; p-1, r)$, where $r$ is an arbitrary integer    \label{fig5:subfig2}}]{
\psfrag{z}{$z$}  
\psfrag{w}{$w$}
  \includegraphics[scale=.35]{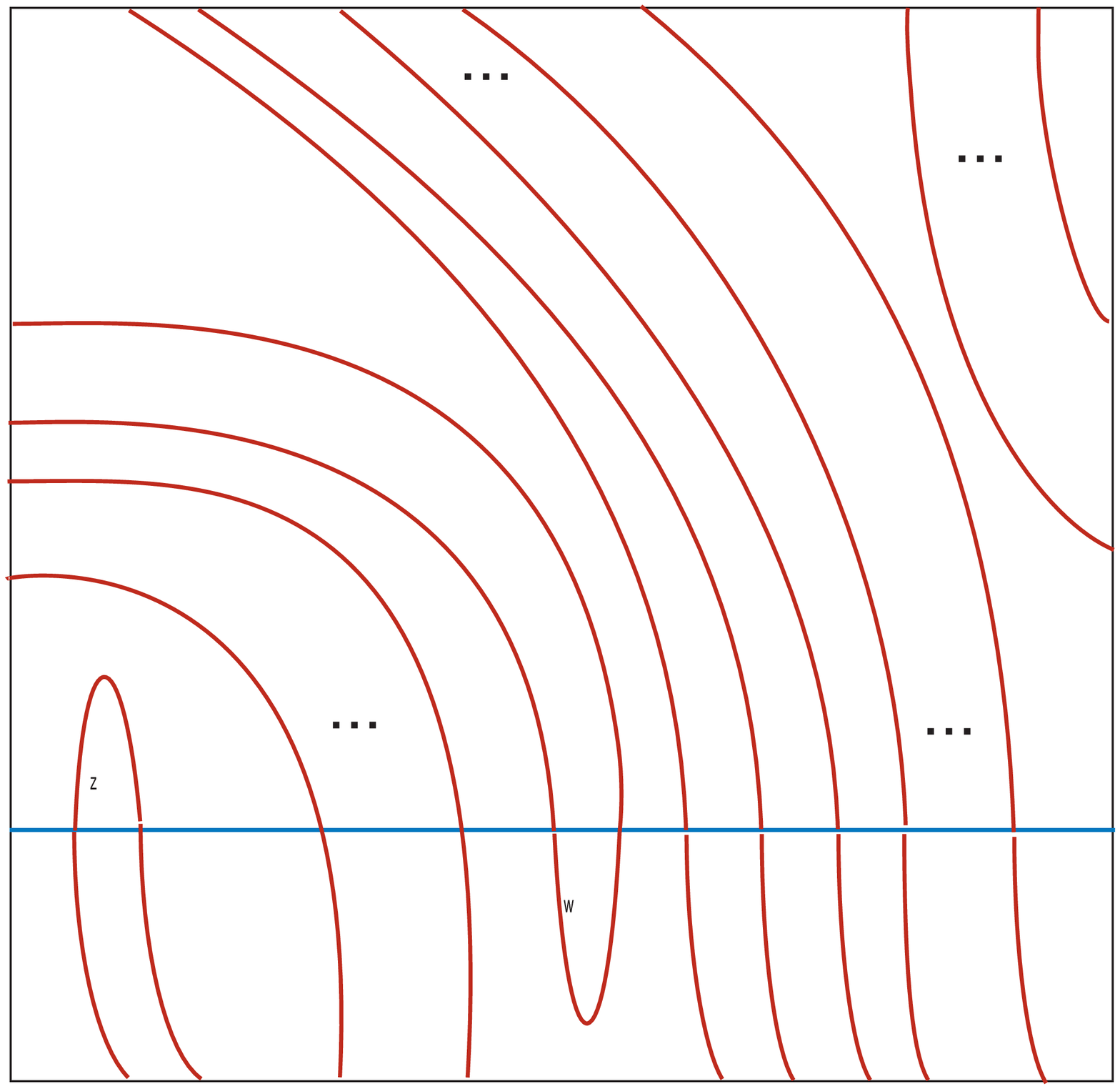}

   }
\caption{\small{Heegaard diagrams on the torus, represented by a rectangle with opposite sides identified}}
       \label{fig5}{%
}
\end{figure}
\\

We prove part (a) by explicitly computing the knot Floer complex of $K(p, k p \pm 1; p-1, r)$. Parts (b) and (c) are proved by focusing on the similarities and differences of the corresponding complexes to those of $K(p, k p \pm 1; p-1, r)$. The key to the proof is in identifying whether or not the knot Floer complex associated to $K(p, k p \pm 1; s, r)$ has a staircase-shape (Corollary \ref{lem1}). 

\begin{proof}[Proof of Theorem \ref{theorem:2}(a)]It will help to break the proof into two steps:

\textbf{Proof Step 1}: We show that $K(p, k p \pm 1; p - 1, r)$ can be presented by a genus one Heegaard diagram with the general form given in  Figure~\ref{fig5:subfig2}. 
 
Case 1: We first consider the case $K(p, k p + 1; p - 1, r)$. The case $p = 2$ is trivial. The construction of a Heegaard diagram in the case when $p = 3$ was given in Section \ref{section:2}. Also Figure \ref{fig14} shows the process for $K = K(4, 5; 3, 2)$. 

To obtain a Heegaard diagram when $p \ge 5$ we can follow a similar procedure. Note that the $w$ base point winds around the longitude of the torus once in the case $p =3$, twice in the case $p=4$, and $p-2$ times in general. Moreover, in each longitudinal winding, the $w$ base point traverses the torus $k + r$ times meridionally, except for the last longitudinal winding where $\alpha$ traverses the torus only $k$ times meridionally. The latter fact holds since we are twisting $p - 1$ strands of the $(p, kp +1)$ torus knot (set $s = p-1$ in Step 3 of the algorithm given in Section \ref{section:2}). Note that as a result of $s = p-1$, we always drag only one sub-arc of $\alpha$ around the torus (Remark \ref{rmk}). Translating the resulting Heegaard diagram obtained this way into the rectangular representation of the torus, we get the general form of Figure~\ref{fig5:subfig2}.

\begin{figure}[t!]
 \centering
 \subfigure[   \label{fig6:subfig1}]{
\psfrag{z}{$z$}
\psfrag{w}{$w$}
  \includegraphics[scale=.4]{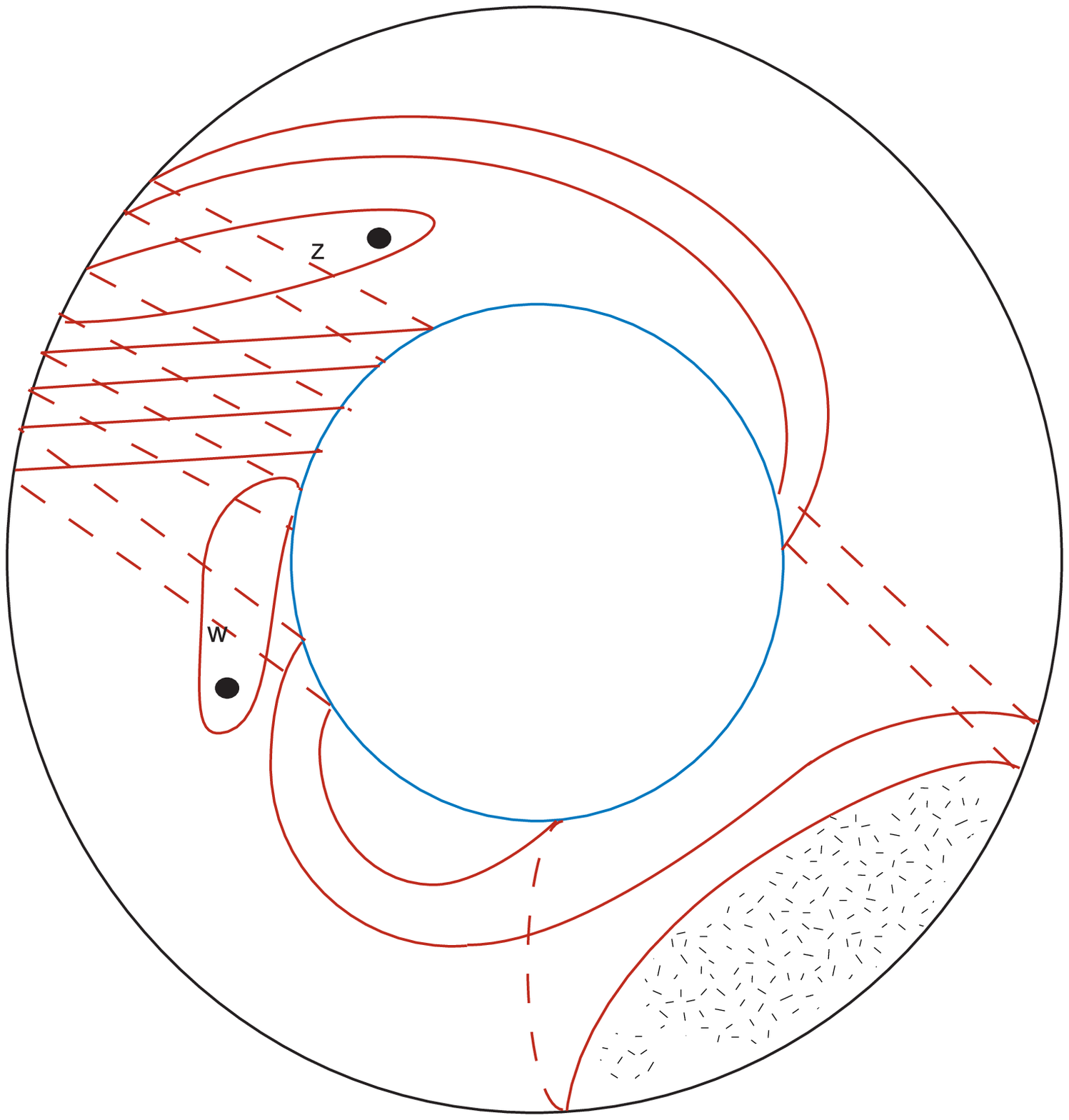}

   }
 \subfigure[   \label{fig6:subfig2}]{
\psfrag{z}{$z$}
\psfrag{w}{$w$}
  \includegraphics[scale=.4]{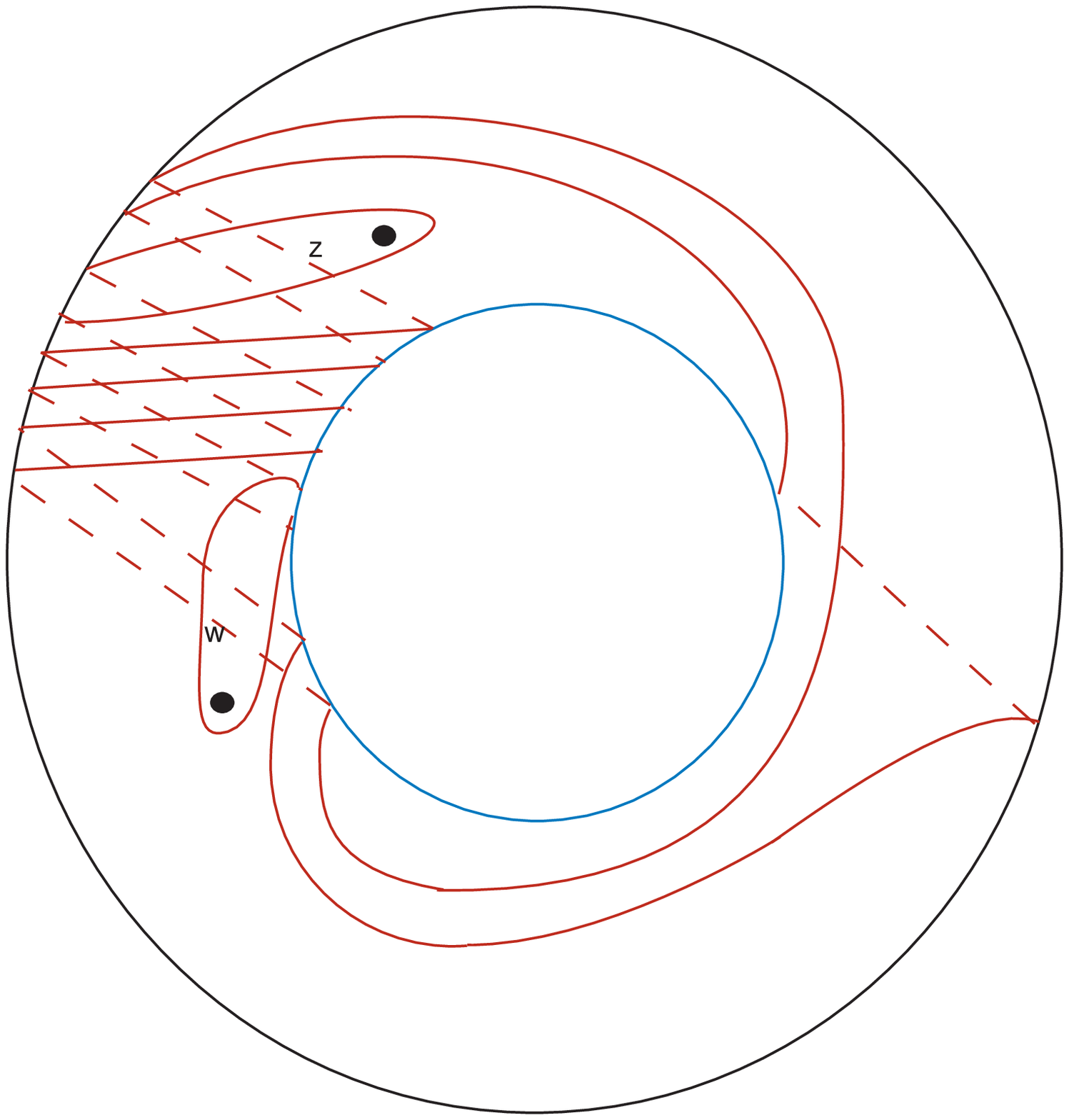}

   }
       \caption{\small{By an isotopy, the shaded region disappears and the Heegaard diagram will have two less intersection points.}} 
\label{fig6}{%
}
\end{figure}

Case 2: For the case $q = kp -1$ we will have a similar setup, though the base points have to pass to the right of each other, not to the left. In this case, there will always be two intersection points of $\alpha$ and $\beta$ that can be removed by an isotopy (see Figure~\ref{fig6:subfig1}). To indicate the general case, we consider $K = K(3, 5; 2, 1)$. The resulting Heegaard diagram is isotopic to a Heegaard diagram for $K(3, 4; 2, 2)$ shown in Figure~\ref{fig6:subfig2}. As in Case 1, the Heegaard diagram will have the general form of Figure~\ref{fig5:subfig2}.

\textbf{Proof Step 2}: In this step, the goal is to calculate the filtered chain complex $CFK^{-}(K)$ for $K = K(p, kp \pm 1; p-1, r)$. Figure \ref{fig8} shows $CFK^{-}(K(3, 4, 2, 2))$. We claim that, in general, $CFK^-(K)$ has the same staircase-shape. 

 As in Section \ref{subsec2} we lift the diagrams, obtained in Step 1, to $\mathbb{C}$ . Fix a connected component $\tilde{\alpha}$ of $\pi^{-1}(\alpha)$. We claim that such a component is a union of $``N"$-shapes (Figure~\ref{fig7:subfig1}). To see this fact, we notice that the lift of a genus one Heegaard diagram can be obtained by gluing together infinitely many copies of the rectangular form of the Heegaard diagram in the plane (gluing from the sides of the rectangles). Figure~\ref{fig7:subfig2} represents a portion of such a lift for a specific example. Pick an intersection point and start moving it along the $\tilde{\alpha}$ curve. (For example, pick the intersection point $9$ on $\tilde{\alpha}$ in Figure~\ref{fig7:subfig2} and start moving it upward.) The direction of the motion will reverse by turning around either of the $z$ or $w$ base points. (In Figure~\ref{fig7:subfig2}, the direction of the motion will change from upward to downward, and also from downward to upward, by going from $1$ to $2$, and from $3$ to $4$, respectively.) Note that the rectangular form of the genus one Heegaard diagram of $K$, as depicted in Figure~\ref{fig5:subfig2}, consists of a single $\beta$ arc, together with $\alpha$ arcs having endpoints on the edge(s) of the rectangle. Note also that there are only two $\alpha$ arcs with both of their endpoints lying on one edge of the rectangle (namely the arcs that turn around the base points). Therefore, by thinking of the lift of the diagram in $\mathbb{C}$ as coming from infinitely many rectangles glued together along the sides and fixing a connected component of $\pi^{-1}(\alpha)$, the change in the direction of the motion (equivalently, turning around either the $z$ or $w$ base point) never happens twice in a single rectangle.\footnote{Note that we do not distinguish between the $z$ and $w$ base points downstairs, and their lifts in $\mathbb{C}$.} Moreover, to recover all the intersection points in the lift, only two changes of direction are needed. As a result, we get the shape of the lifted digram as claimed.  

Let us first consider the example, $CFK^{-}(K(3, 4; 2, 2))$ whose Heegaard diagram is given in Figure~\ref{fig7:subfig2}. Given a pair of intersection points $x$ and $y$, the moduli space of holomorphic representatives of Whitney disks $\phi \in \pi_2(x, y)$ with Maslov index one, modulo reparametrization, is either empty or consists of one map. In what follows, we write $x \rightarrow y$ if the moduli space consists of one such map, and if so, we record how many times it passes over the $z$ and $w$ base points:
\begin{itemize}

\item $2 \rightarrow 1$, $6 \rightarrow 5$, $8 \rightarrow 7$   using one $z$ base point,

\item  $3 \rightarrow 9$ using two $z$ base points, 

\item $6 \rightarrow 7$, $8 \rightarrow 9$, $3 \rightarrow 4$  using one $w$ base point, and

\item $2 \rightarrow 5$ using two $w$ base points.
\end{itemize}

From Figure~\ref{fig7:subfig2}, it is easy to see that we need four $\tilde{\beta}$ lines to generate the whole nine intersection points in the lifted Heegaard diagram, i.e. fixing $\tilde{\alpha}$, by using only four connected components of the lift of $\beta$ we can obtain a lift of all the intersection points between $\alpha$ and $\beta$. Starting from $\tilde{\beta}_4$,

\newpage
\begin{figure}[htb!]
 \centering
 \subfigure[\label{fig7:subfig1}]{
\psfrag{z}{$z$}
\psfrag{w}{$w$}
\psfrag{1}{\tiny $x_{m-2}$}
\psfrag{2}{\tiny $x_{m-1}$}
\psfrag{3}{\tiny $x_m$}
\psfrag{4}{\tiny $x_{m-2}$}
\psfrag{5}{\tiny $x_{m-3}$}
\psfrag{6}{\tiny $x_{m-4}$}
\psfrag{7}{\tiny $x_3$}
\psfrag{8}{\tiny $x_2$}
\psfrag{9}{\tiny $x_1$} 
\psfrag{a}{\tiny $\tilde{\beta}_n$}
\psfrag{b}{\tiny $\tilde{\beta}_{n-1}$}
\psfrag{c}{\tiny $\tilde{\beta}_2$}
\psfrag{d}{\tiny $\tilde{\beta}_1$}
\includegraphics[scale=.37]{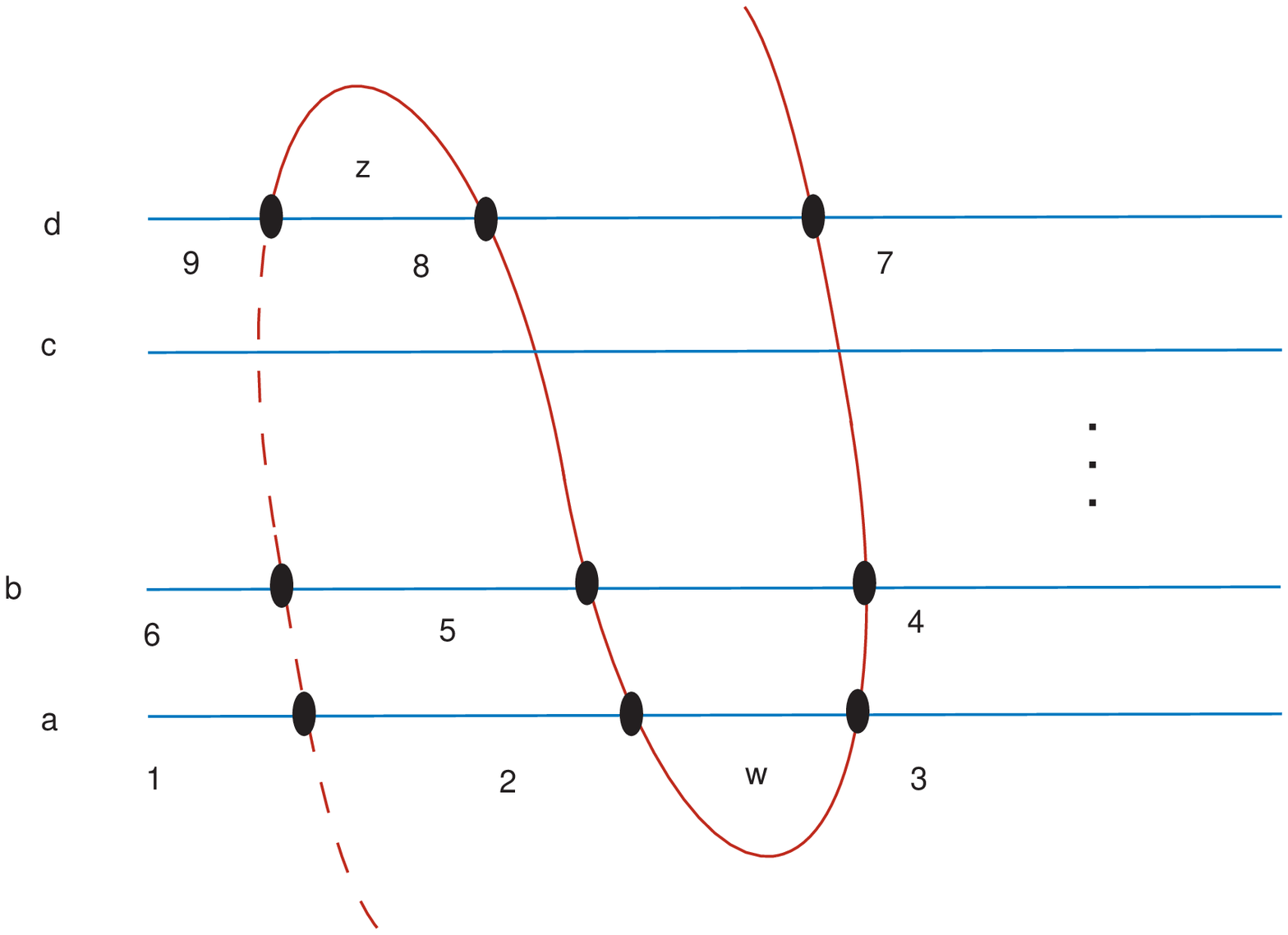}

   }
 \subfigure[\label{fig7:subfig2}]{
\psfrag{z}{$z$}  
\psfrag{w}{$w$}  
\psfrag{1}{\tiny$1$}
\psfrag{2}{\tiny$2$}
\psfrag{3}{\tiny$3$}
\psfrag{4}{\tiny$4$}
\psfrag{5}{\tiny$5$}
\psfrag{6}{\tiny$6$}
\psfrag{7}{\tiny$7$}
\psfrag{8}{\tiny$8$}
\psfrag{9}{\tiny$9$} 
\psfrag{a}{$\tilde{\beta}_1$}
\psfrag{b}{$\tilde{\beta}_2$}
\psfrag{c}{$\tilde{\beta}_3$}
\psfrag{d}{$\tilde{\beta}_4$}
  \includegraphics[scale=.28]{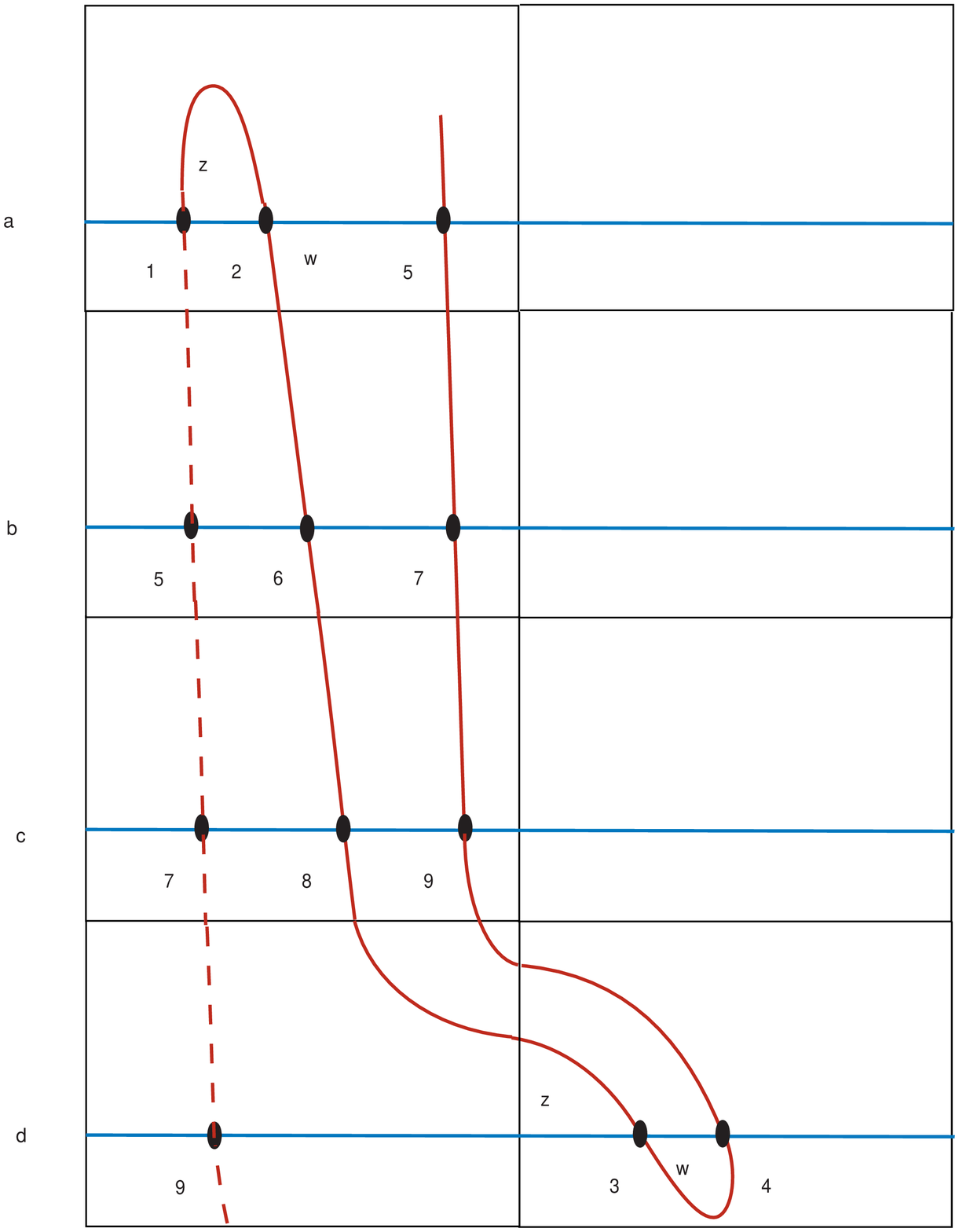}

   } 
\caption{\small{(a) A portion of the Heegaard diagram for $K=K(p, kp \pm 1; p-1, r)$ lifted to $\mathbb{C}$, where $r$ is an arbitrary integer. Note that $m$ is the number of intersection points in the genus one Heegaard diagram of $K$. It is assumed, fixing $\tilde{\alpha}$ a connected component of $\pi^{-1}(\alpha)$, that we need $n$ connected components of $\pi^{-1}(\beta)$ to obtain a complete list of all the $m$ intersection points between $\alpha$ and $\beta$ downstairs. (b) A portion of the Heegaard diagram for the $(3, 4)$ torus knot with two positive full twists on two adjacent strands, lifted to $\mathbb{C}$. Note that the base points specified in the picture depicted above are the only relevant base points needed to compute $CFK^-.$}}
\label{fig7}{
}
\end{figure}
\begin{figure}[H]
\begin{center}
\psfrag{2}{\tiny$1$}
\psfrag{1}{\tiny$2$}
\psfrag{4}{\tiny$3$}
\psfrag{3}{\tiny$4$}
\psfrag{5}{\tiny$5$}
\psfrag{6}{\tiny$6$}
\psfrag{7}{\tiny$7$}
\psfrag{8}{\tiny$8$}
\psfrag{9}{\tiny$9$}   
\includegraphics[ scale=.241]{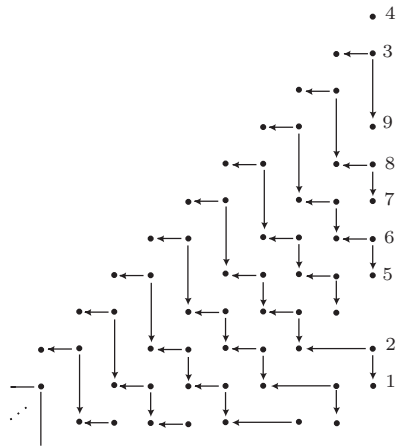}
\end{center}
\caption{\small{$CFK^-(K(3, 4, 2, 2))$}}
\label{fig8}
\end{figure}
\newpage
\noindent 
 there are three intersection points (3, 4 and 9) with one disk $4 \rightarrow 3$ using one $w$ base point and one other disk $9 \rightarrow 3$ using two $z$ base points. Thus, in terms of the Alexander gradings $A(i)$ of the intersection points, $i \in \left\{ 1, 2, ..., 9\right\}$ , we have that:
\begin{itemize}

\item $A(3) - A(4) = n_{z}(\tilde{\phi}) - n_{w}(\tilde{\phi}) = -1$, and

\item $A(3) - A(9) = n_{z}(\tilde{\phi}) - n_{w}(\tilde{\phi}) = 2$.
\end{itemize}
\noindent
See \cite{Ozsvath2004} for the notation. By a similar method, we can find the remaining Whitney disks in the 
 list above and use them to complete the ordering of the Alexander gradings. At this point, we can obtain the staircase-shape of Figure \ref{fig8}. 

For the general case of Figure~\ref{fig7:subfig1}, it is straightforward to observe that our strategy can be extended. Assume that $\{ x_1, x_2, ..., x_m \}$ is the set of intersection points between $\alpha$ and $\beta$ curves coming from the genus one Heegaard diagram of $K = K(p, kp+1; p-1, r)$ (see Figure~\ref{fig5:subfig2} and Figure~\ref{fig7:subfig1}). Assume also that, fixing $\tilde{\alpha}$ a connected component of $\pi^{-1}(\alpha)$, we need $n$ connected components of $\pi^{-1}(\beta)$ to recover all the $m$ intersection points downstairs between $\alpha$ and $\beta$ (Figure~\ref{fig7:subfig1}). Our strategy is first ordering the generators based on their Alexander gradings and, second, finding all the differentials. Using the ``N"-shape of Figure~\ref{fig7:subfig1} and starting from $\tilde{\beta}_n$, there are three intersection points ($x_m$, $x_{m-1}$ and $x_{m-2}$) with one disk $x_{m-1} \rightarrow x_{m}$ using one $w$ base point and one other disk $x_{m-1} \rightarrow x_{m-2}$ using the $z$ base point(s). Note that there exists no other non-trivial Whitney disk with Maslov index one connecting $x_{m-1}$ to another intersection point of Figure~\ref{fig7:subfig1}. Also on $\tilde{\beta}_{n-1}$, there is one disk $x_{m-3} \rightarrow x_{m-2}$ using the $w$ base point(s). Continuing this process, we deduce that 
\[ A(x_m) > A(x_{m-1}) > A(x_{m-2}) >  A(x_{m-3}) > ... >A(x_1). \] 
By noting that there is no other non-trivial Whitney disk with Maslov index one, we see that the set $\{ x_1, x_2, ..., x_m\}$ forms a basis for $CFK^{-}(K)$ such that 
\[
\begin{array}{ccc}
\partial x_i = &x_{i-1}+x_{i+1} & \text{for $i$ even} \\
\partial x_i = &0                      & \text{otherwise.}
\end{array}
\]
\noindent This formula for the differentials (which is the same as \eqref{staircase}), together with the existence of three intersection points on each $\tilde \beta_{j}$ line of Figure~\ref{fig7:subfig1} with exactly two disks using different base point types (i.e. $z$ and $w$), gives the staircase-shape of $CFK^-(K)$ (see the discussion about a knot Floer complex with a staircase-shape in Section~\ref{subsec2}). Now, Corollary \ref{lem1} completes the proof.

\end{proof}
\noindent
\begin{proof}[Proof of (b) and (c)]Let $K(p, q; s, r)$ be a twisted torus knot where $2 \le s \le p-2$. We discuss the case when $q = kp+1$ and leave the case $q=kp-1$ to the reader. Since we apply the same algorithm, as used in Part (a), to obtain a Heegaard diagram, we will only highlight the differences in this case. Recalling the algorithm explained in Section \ref{subsec1}, we first wind $z$ once in the counterclockwise direction (Step 1). Then we wind the $w$ base point $(s-2)$ times in the clockwise direction, traversing the torus $(k + r)$ times meridionally in each winding (Step 2). Finally, we wind the $w$ base point $(p - s)$ more times around the torus longitudinally (Step 3). Note that in the latter step, $w$ goes through only $k$ meridional moves in each winding. 

 It will be convenient to pick an arbitrary orientation for the $\alpha$ curve. Note that, unlike Part (a), more than one sub-arc will be dragged since $2 \le s \le p-2$ (Remark \ref{rmk}). With the $\alpha$ curve oriented, either these sub-arcs will have all the same orientation or there will be at least one pair of sub-arcs with opposite orientations. The case for only two sub-arcs can be seen in Figure~\ref{fig120}. Figure \ref{fig18} shows the process of constructing a Heegaard diagram for $K(4, 5, 2, 1)$, which indicates the pattern, particularly in the case when $s \in \left\{ 2, p-2 \right\}$.

\newpage

\begin{figure}[H]
 \centering
 \subfigure[   \label{fig14:subfig1}]{
  \psfrag{z}{$z$}  
\psfrag{w}{$w$}
  \includegraphics[scale=.3]{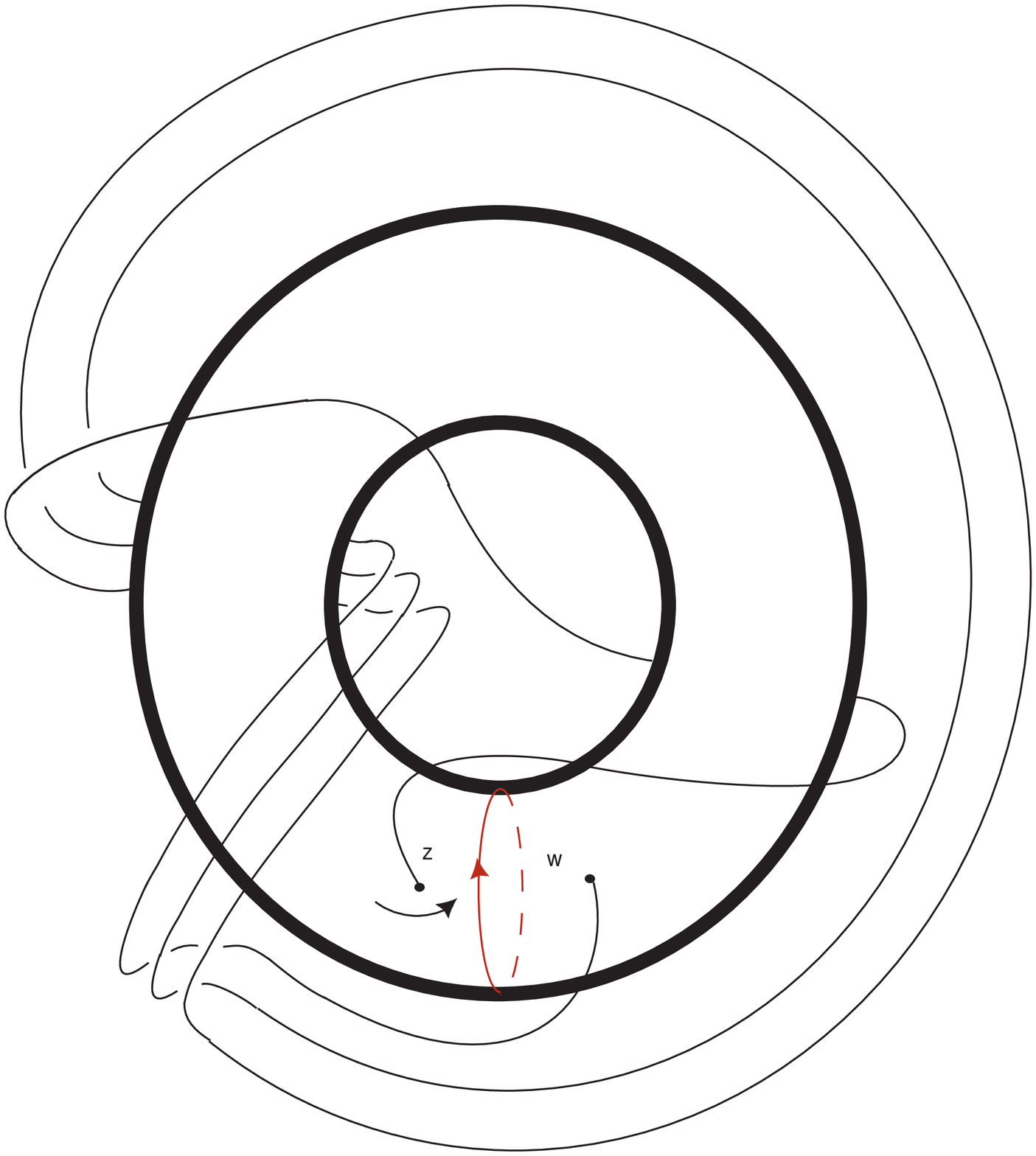}

   } 
\subfigure[   \label{fig14:subfig2}]{
  \psfrag{z}{$z$}  
\psfrag{w}{$w$}
  \includegraphics[scale=.3]{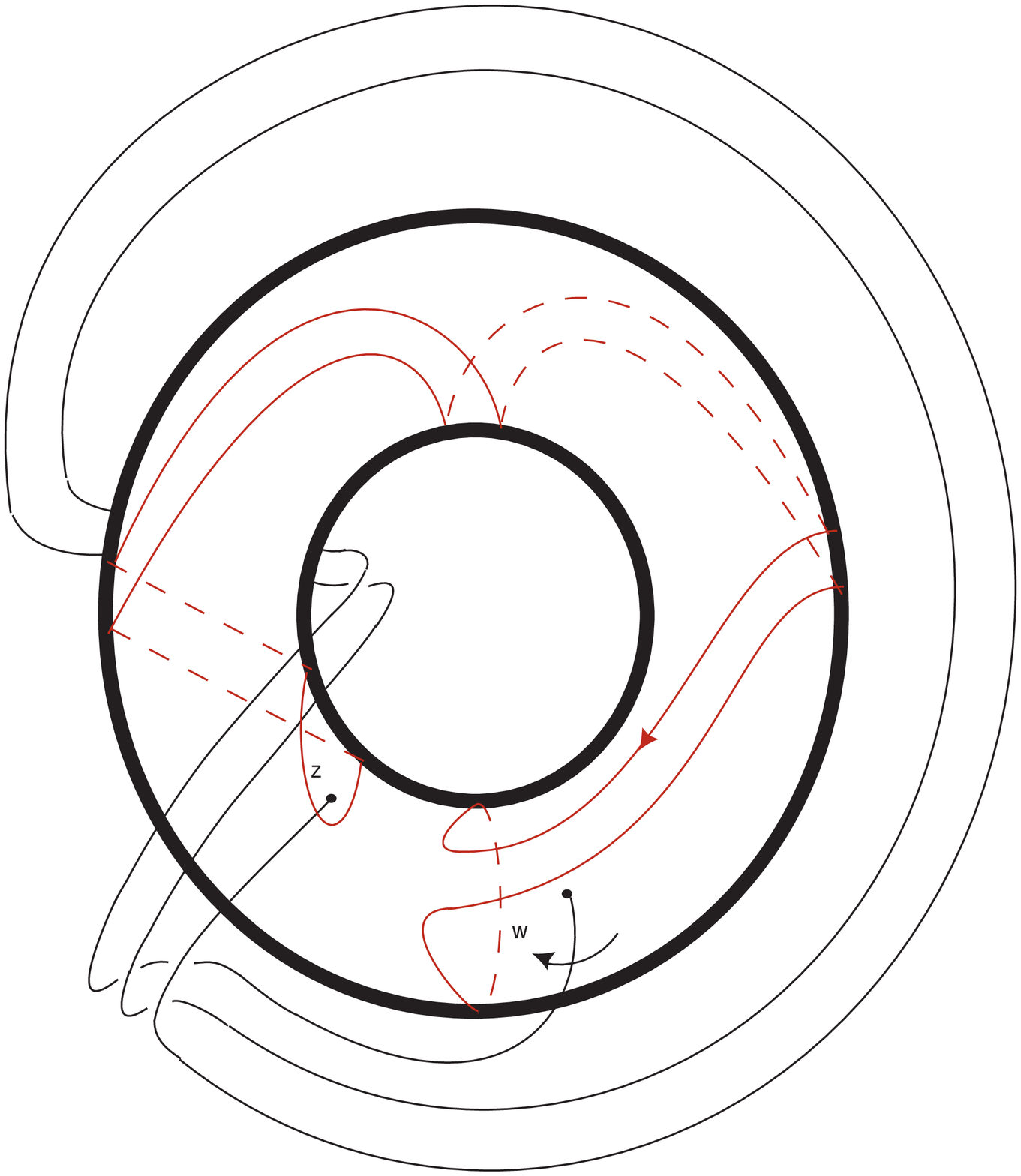}

   }
 \subfigure[   \label{fig14:subfig3}]{
  \psfrag{z}{$z$}  
\psfrag{w}{$w$}
  \includegraphics[scale=.3]{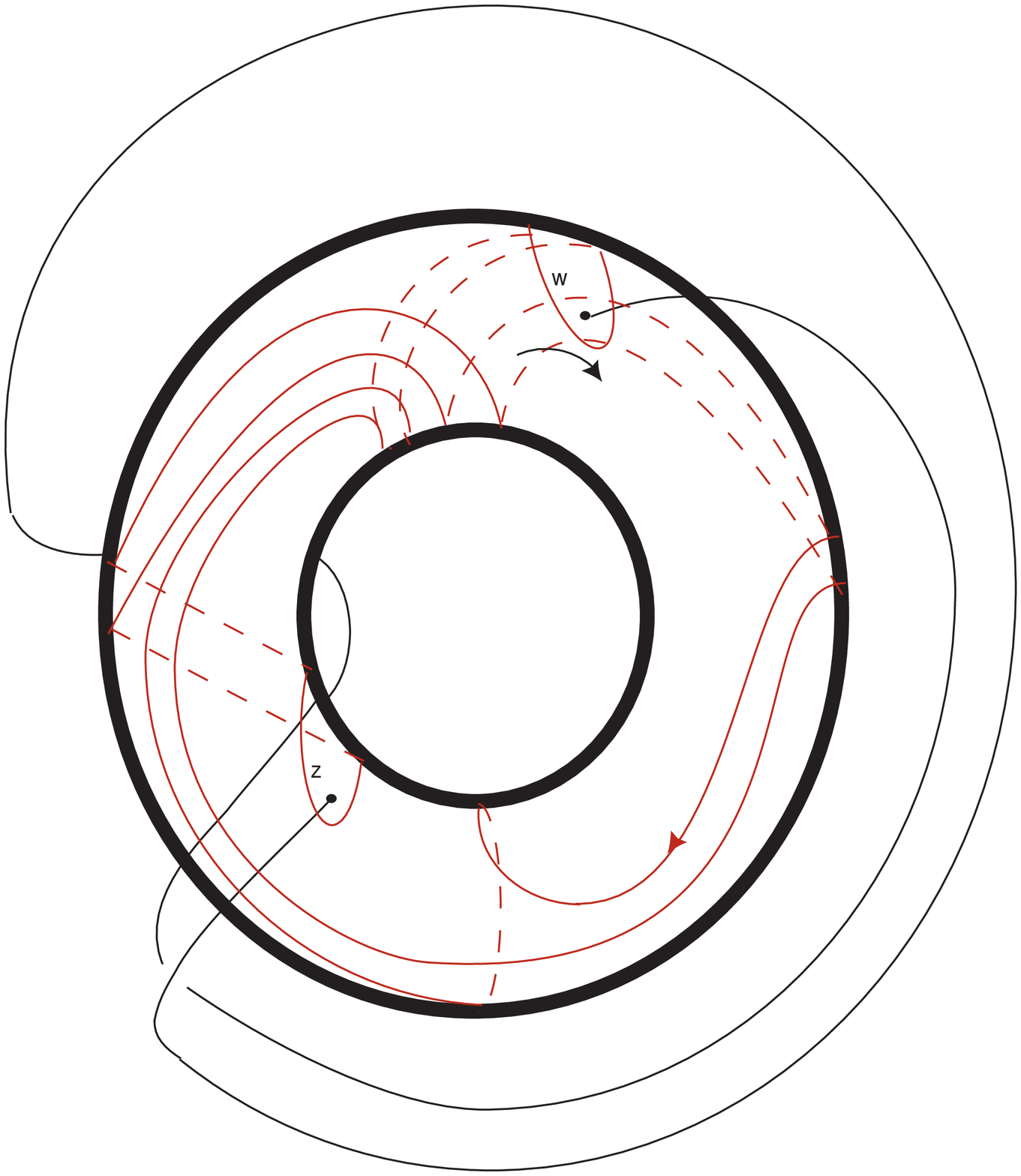}

   }
 \subfigure[   \label{fig14:subfig4}]{
  \psfrag{z}{$z$}  
\psfrag{w}{$w$}
  \includegraphics[scale=.3]{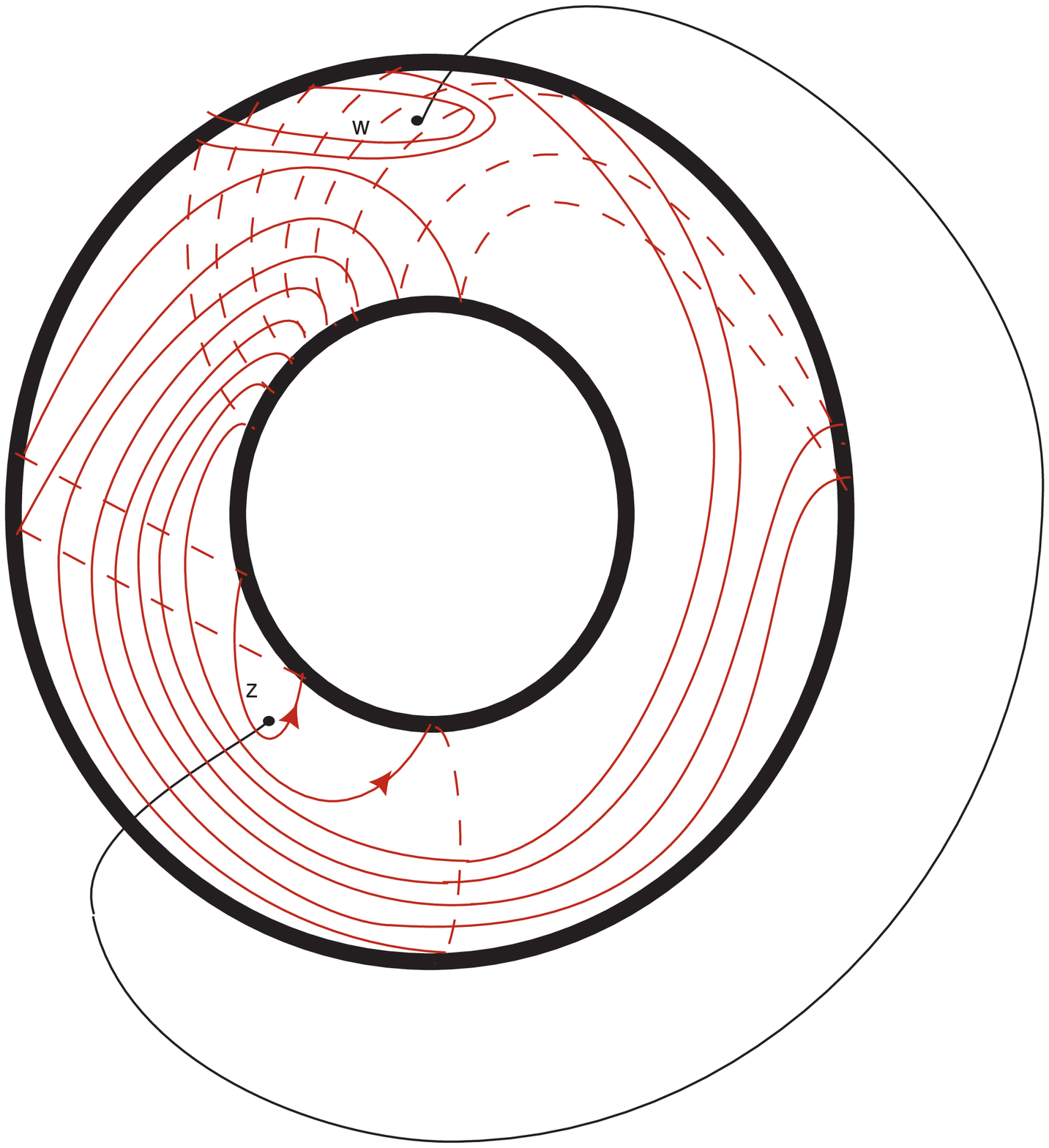}

   }

\caption{\small{The process of drawing a genus one Heegaard diagram for $K(4, 5; 2, 1)$. The $\alpha$ curve in each step is oriented. This example indicates the pattern when $s \in \left\{ 2, p-2 \right\}$ and $r = 1$. In general when $r=1$, to go from (c) to (d), $w$ first drags $(s - 1)$ sub-arcs, all oriented in the same direction. In the next winding it drags $(s - 2)$ additional sub-arcs, all oriented in the same direction but opposite to those of the first $(s-1)$ sub-arcs. Dragging oppositely oriented sub-arcs does not occur in this example since $s = 2$. Note that the orientation is irrelevant once the Heegaard diagram is completed.}}{%
 }
\label{fig18}{%
}
\end{figure}
\newpage
\begin{figure}[H]
\centering
\psfrag{w}{$w$}
\includegraphics[scale = .4]{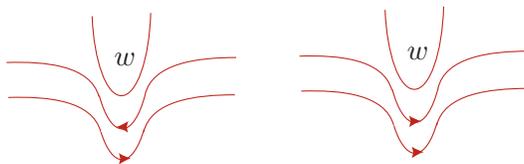}
\caption{\small{The base point $w$ drags more than one sub-arc of $\alpha$. The picture depicted above is schematic.}}
\label{fig120}
\end{figure}



\noindent \textbf{Claim.} Unless $s\in \left\{2, p-2\right\}$ and $r = 1$, the trivializing process will drag oppositely oriented sub-arcs.

\begin{proof}
\noindent Suppose $r = 1$. The first longitudinal traversal of Step 3 drags no additional sub-arcs. The second traversal of Step 3, however, drags $(s - 1)$ sub-arcs, all oriented in the same direction. The next winding drags $(s - 2)$ additional sub-arcs, all oriented in the same direction but opposite to those of the first $(s-1)$ sub-arcs. This opposite orientation will clearly not occur if $s = 2$. Suppose $s = p-2$. Then in Step 3 the $w$ base point is wound longitudinally around the torus $p - (p-2) = 2$  times (twice). Hence, only sub-arcs with the same direction will be dragged.  If $r \ge 2$ the full twists of Step 1 create future oppositely oriented sub-arcs in Step 3, i.e. the $w$ base point will be dragging sub-arcs with opposite orientations, starting the second longitudinal traverse of Step 3.  More specifically, if the number of full twists is greater than one, each additional twist will create two oppositely oriented sub-arcs and the $w$ base point will drag both of these sub-arcs after the first $(s-1)$ longitudinal windings.
\end{proof}
Since the hypotheses of Part (b) imply that the sub-arcs have the same orientation, a similar argument to Part (a), once we lift the diagram to $\mathbb{C}$, shows that the ordering of the Alexander gradings of the intersection points will follow the same manner as in the case $s = p-1$. More precisely, if we think of the lift of the Heegaard diagram as coming from infinitely many rectangles glued together, by picking an intersection point and moving it along a fixed connected component $\tilde{\alpha}$ of $\pi^{-1}(\alpha)$, we see that the picked point, during its motion, never turns around the $z$ (or $w$) base point twice in a single rectangle. Therefore, although the lifted diagrams are not looking the same as Part (a), we claim that the corresponding complexes have the staircase-shape. In particular, for the case $s = 2$ (respectively $s = p-2$), we need four (respectively $2p-4$) changes of direction\footnote{We remind the reader that by changing direction we mean turning around one of the base points ($z$ or $w$).} to recover all the intersection points of downstairs. For the specific example of $K(4, 5; 2, 1)$ depicted in Figure~\ref{fig111:subfig2}:
\[A(6)>A(5)>A(9)>A(8)>A(7)>A(4)>A(1)>A(11)>A(10)>A(3)>A(2).\]
Exploring the Whitney disks in the lifted diagram will give a staircase-shape for the associated complex. To see this in the general case, note that the set of all intersection points between $\alpha$ and $\beta$ curves forms a basis for $CFK^-$. Moreover, for every intersection point $x_i$, either the differential vanishes, or there exist two Whitney disks with Maslov index one connecting $x_i$ to another two distinct intersection points, using $z$ and $w$ base points alternatively. (This shows that the differentials are of the form of \eqref{staircase}.) That is, for each intersection point $x_i$, either there is no arrow joining it to another intersection point, or there are two arrows joining $x_i$ to two distinct intersection points such that one arrow is horizontal and the other is vertical. This gives us the staircase-shape of the knot Floer complex. Finally, Corollary \ref{lem1} completes the proof of Part (b).
\\

\noindent 
To prove Part (c), note that if the arcs dragged by $w$ have different orientations, then, after lifting the diagram to $\mathbb{C}$, the following phenomenon occurs:
\begin{figure}[H]
\centering
\psfrag{w}{$w$}
\psfrag{z}{$z$}
\psfrag{1}{\tiny{$1$}}
\psfrag{2}{\tiny{$2$}}
\psfrag{3}{\tiny{$3$}}
\psfrag{4}{\tiny{$4$}}
\psfrag{a}{\tiny{$\tilde{\beta}_i$}}
\includegraphics[scale = .3]{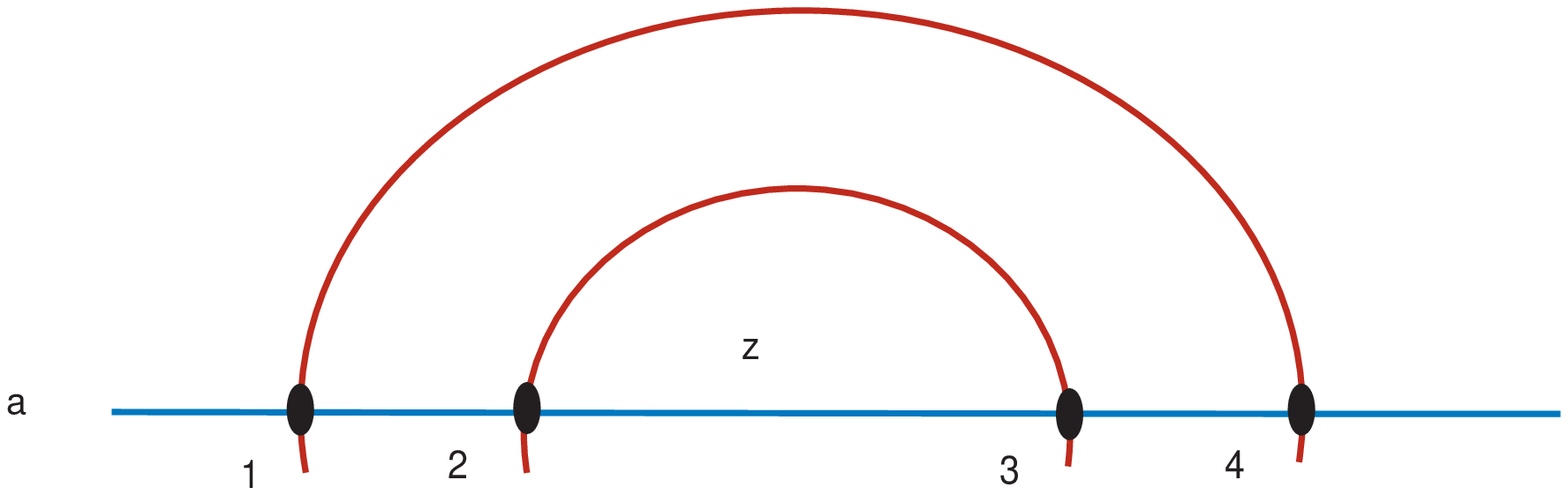}
\end{figure}
\noindent \textbf{Claim:} The associated complex does not represent an L-space knot.
\\

\noindent \emph{Proof.} As in the proof of Part (a), we can order the Alexander gradings of the intersection points from the Whitney disks in the lifted Heegaard diagram. Let $\tilde{\beta}_1$, ..., $\tilde{\beta}_k$ denote the lifts of $\beta$ needed to find all of the Whitney disks. Work from $\tilde{\beta}_k$ to $\tilde{\beta}_1$ and stop at the first $\tilde{\beta}_i$ that exhibits the phenomenon in Figure above. Then part of the diagram is as Figure~\ref{fig112}.
\begin{figure}[H]
\centering
\psfrag{w}{$w$}
\psfrag{z}{$z$}
\psfrag{1}{\tiny{$1$}}
\psfrag{2}{\tiny{$2$}}
\psfrag{3}{\tiny{$3$}}
\psfrag{4}{\tiny{$4$}}
\psfrag{a}{\tiny{$\tilde{\beta}_k$}}
\psfrag{b}{\tiny{$\tilde{\beta}_{i+1}$}}
\psfrag{c}{\tiny{$\tilde{\beta}_i$}}

\includegraphics[scale = .34]{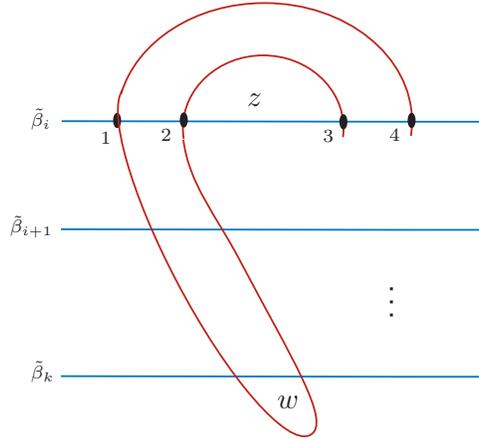}
\caption{\small{A sub-diagram of a lifted Heegaard diagram, fixing one connected component of $\tilde{\alpha}$ }}
\label{fig112}
\end{figure}
\noindent We analyze this by looking at the Whitney disks:
\begin{itemize}

\item $4 \rightarrow 1$, $3 \rightarrow 2$ using one $z$ base point, and

\item $1 \rightarrow 2$, $4 \rightarrow 3$ using one $w$ base point.

\end{itemize} 
As a result, the part of $CFK^-$ involving the intersection points, $\left\{1, 2, 3, 4\right\}$, on $\tilde{\beta}_i$ will look like
\begin{figure}[H]
\centering
\psfrag{w}{$w$}
\psfrag{z}{$z$}
\psfrag{1}{\tiny{$1$}}
\psfrag{2}{\tiny{$2$}}
\psfrag{3}{\tiny{$3$}}
\psfrag{4}{\tiny{$4$}}
\includegraphics[scale = .38]{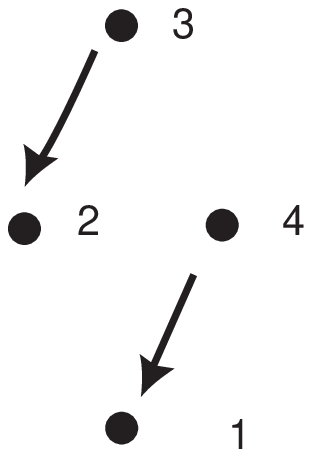}

\end{figure}

\noindent Note that the boundary map decreases the Maslov grading by one, and the $U$-action decreases the grading by two. Combining these facts with the existence of the disks $1 \rightarrow 2$ and $4 \rightarrow 3$, we find that the intersection points 2 and 4 both have the same Maslov gradings as well as the same Alexander gradings. (We are assuming that there are no trivial Whitney disks connecting two intersection points; if there is a bigon that does not pass over any of the base points, we can isotop it away.) Thus,
\[
\rank \text{ }\widehat{HFK}(K, \mathfrak{s}) \ge \rank \text{ }\widehat{HFK}_m(K, \mathfrak{s}) = \rank \text{ }\widehat{CFK}_m(K, \mathfrak{s}) \ge 2,
\]
\noindent where $\mathfrak{s}$ is the Alexander grading of the intersection points $2$ and $4$. Now, Lemma~\ref{lem0} completes the proof of the claim and Part (c).
\end{proof}
The Heegaard diagrammatic observation in Figure \ref{fig6} can be generalized. The author suspects that the following corollary could have been proved differently, using braid words for instance: 
\begin{figure}[t!]
\subfigure[\label{fig111:subfig1}]{
\psfrag{1}{\tiny{$1$}}
\psfrag{2}{\tiny{$2$}}
\psfrag{3}{\tiny{$3$}}
\psfrag{4}{\tiny{$4$}}
\psfrag{5}{\tiny{$5$}}
\psfrag{6}{\tiny{$6$}}
\psfrag{7}{\tiny{$7$}}
\psfrag{8}{\tiny{$8$}}
\psfrag{9}{\tiny{$9$}}   
\psfrag{10}{\tiny{$10$}}
\psfrag{11}{\tiny{$11$}}
  \psfrag{z}{\tiny{$z$}}  
\psfrag{w}{\tiny{$w$}}  

 \includegraphics[scale=.26]{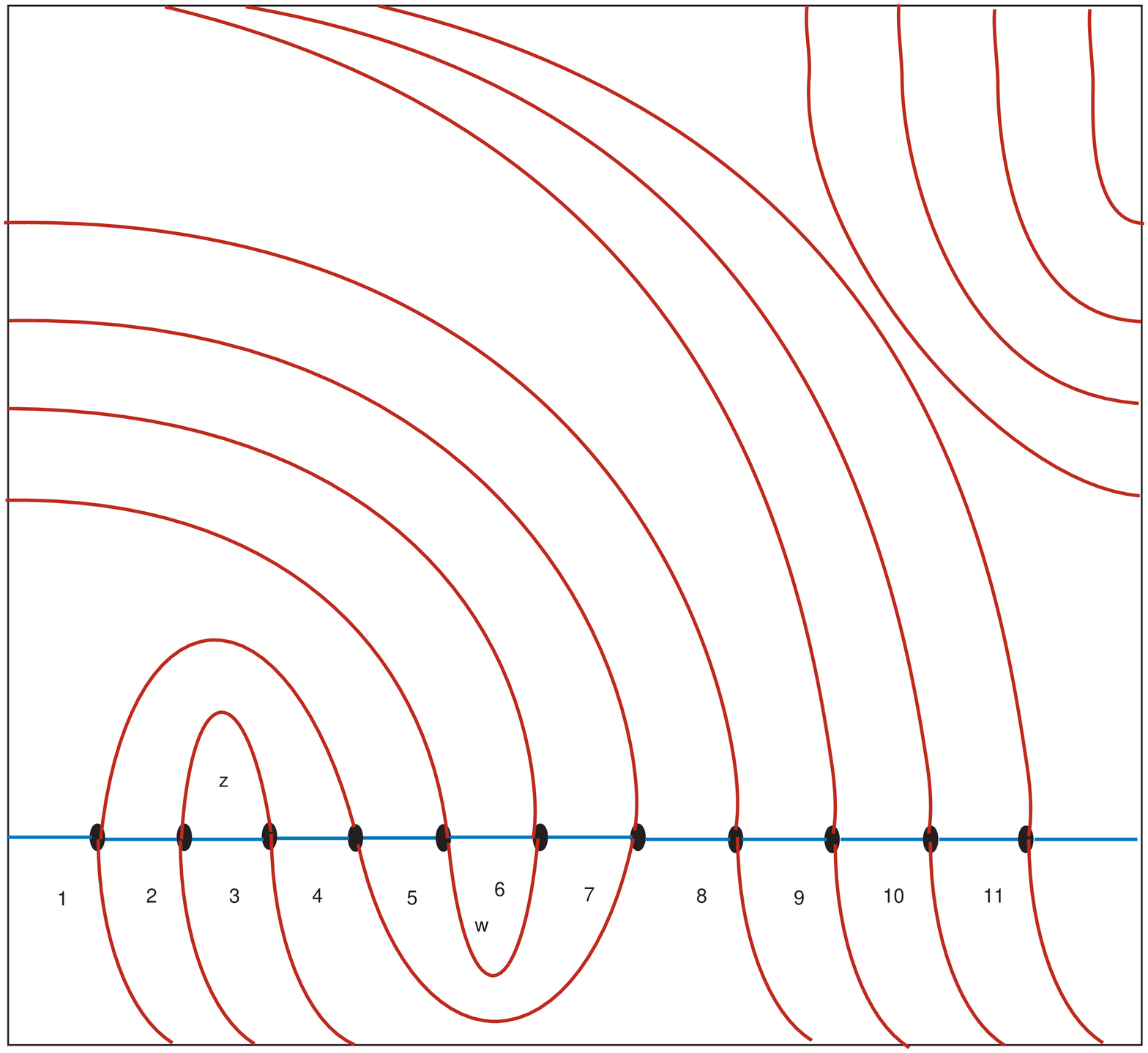}

}
\subfigure[\label{fig111:subfig2}]{
\psfrag{1}{\tiny{$1$}}
\psfrag{2}{\tiny{$2$}}
\psfrag{3}{\tiny{$3$}}
\psfrag{4}{\tiny{$4$}}
\psfrag{5}{\tiny{$5$}}
\psfrag{6}{\tiny{$6$}}
\psfrag{7}{\tiny{$7$}}
\psfrag{8}{\tiny{$8$}}
\psfrag{9}{\tiny{$9$}}   
\psfrag{10}{\tiny{$10$}}
\psfrag{11}{\tiny{$11$}}
  \psfrag{z}{\tiny{$z$}}  
\psfrag{w}{\tiny{$w$}}
\psfrag{a}{\tiny{$\tilde{\beta}_1$}}
\psfrag{b}{\tiny{$\tilde{\beta}_{2}$}}
\psfrag{c}{\tiny{$\tilde{\beta}_3$}}   
\psfrag{d}{\tiny{$\tilde{\beta}_4$}} 
 \includegraphics[scale=.22]{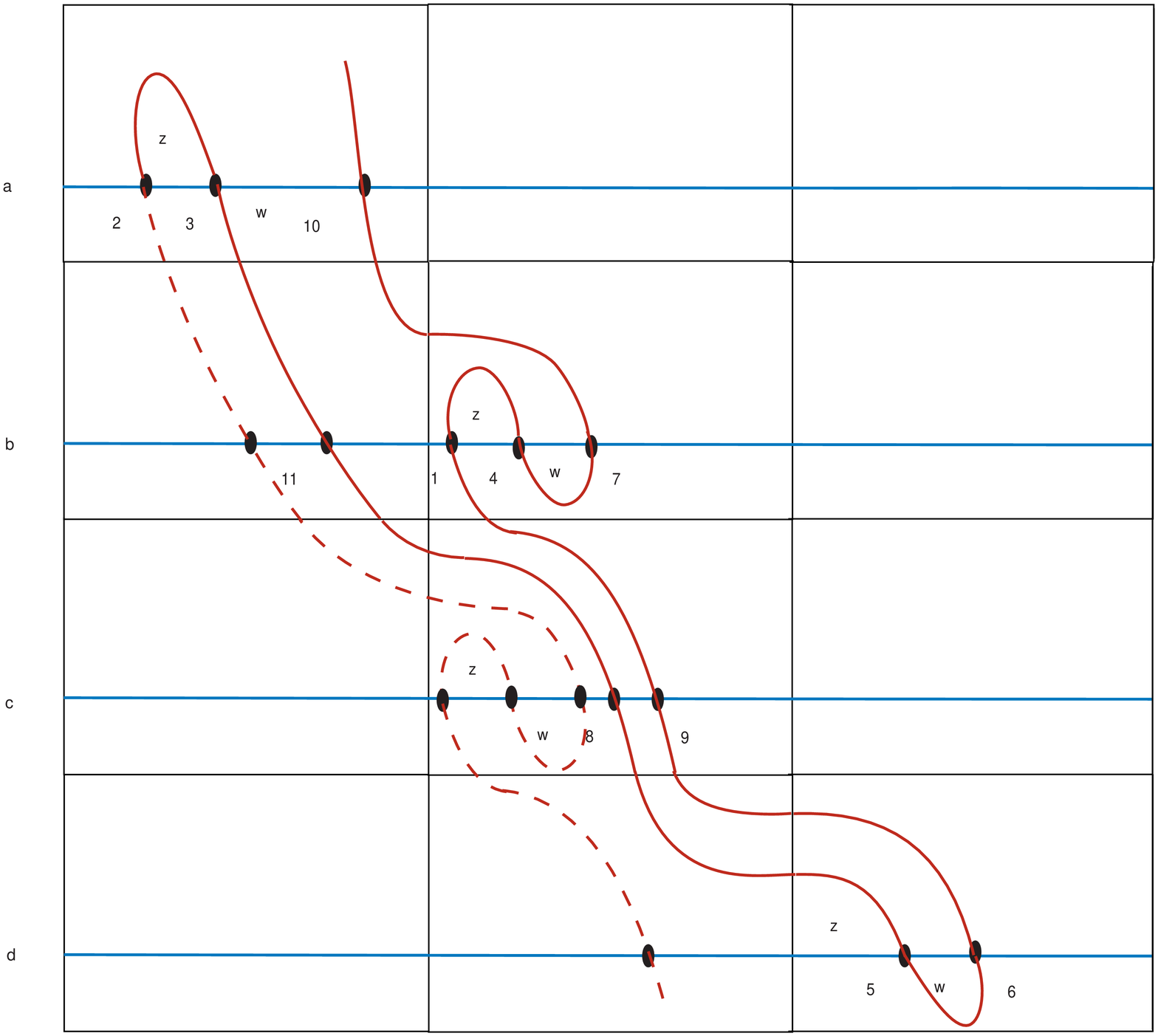}

}
       \caption{\small{A genus one Heegaard diagram for $K(4, 5; 2, 1)$, as well as its lift to $\mathbb{C}$}}
\label{fig111}{%
}
\end{figure}
{\cor \label{cor2}The twisted torus knot, $K(p, kp + 1; p-1, r)$, is isotopic to $K(p, (k + 1)p -1; p-1, r-1)$.}
\begin{proof}
We start from the genus one Heegaard diagram of $K_1 = K(p, (k+1)p + 1; p-1, r-1)$, obtained from implementing the algorithm explained in Section~\ref{subsec2}. The proof is done by first doing an isotopy to get rid of the two extra generators in the genus one Heegaard diagram of $K_1$\footnote{Note that the phenomenon (of having two removable intersection points) in Figure~\ref{fig6}, once we implement the algorithm explained in Section~\ref{subsec2}, will always occur in the genus one Heegaard diagram of $K(p, kp - 1; p-1, r)$.} and, second, tracking back the drag of the $w$ and $z$ base points in the torus. More precisely, after removing the extra generators, if we track back the $w$ base point, we see that it passes, during its $p-2$ longitudinal windings, to the right of $z$. Now, by tracking back the $z$ base point once around the torus, we see that it also passes to the right of $w$. These facts can be verified in the example depicted in Figure~\ref{fig6:subfig2}. (Thus, while implementing the algorithm to obtain the diagram in the first place, the base points must have passed by the left of each other). During this process, except for the first winding of $w$ that goes through $k$ meridional moves, the rest of windings traverse the torus $k+r$ times meridionally. Therefore, by noting that only one sub-arc of $\alpha$ has been dragged by the base points, we get that the diagram obtained after doing the isotopy is a genus one Heegaard diagram for $K_2 = K(p, kp + 1; p-1, r)$.   
\end{proof}
\noindent When $p = 3$ in Theorem~\ref{theorem:1}, we obtain a generalization of \cite[Theorem 1.2]{Watson}: 
{\cor \label{cor3}All twisted $(3, q)$ torus knots admit L-space surgeries.}

\section{Directions for future research} \label{section:4}
Closely related to the main result of this paper, one can ask the question of which operations on knots produce L-space knots. Satellite operations are the first in line. As pointed out in Section~\ref{section:1}, the $(p, q)$ cabling is an L-space satellite operation \cite{Hom2011a}. More generally, Hom, Lidman and the author introduced an L-space satellite operation, using \emph{Berge-Gabai knots} as the pattern \cite{Homa}. By definition, a knot $P \subset S^1 \times D^2$ is called a Berge-Gabai knot if it admits a non-trivial solid torus surgery. We also suspect that one can obtain more L-space satellite operations, choosing the patterns from the the list of L-space knots of Theorem \ref{theorem:1}. Although classifying such operations does not seem to be an easy task to do, there is an obstruction to obtaining L-space satellite knots (Lemma~\ref{lem0}) which can be appealing. Let $P(K)$ be a satellite knot with pattern $P \subset V = S^1\times D^2$ and companion $K$. We recall the behavior of the Alexander polynomial of a satellite knot: 
\[
\Delta_{P(K)}(T) = \Delta_{P}(T) \Delta_{K}(T^w)
\]
where $w$ is the geometric intersection number of the pattern $P$ with a fixed meridional disk of $V$ (see for instance \cite{Lickorish1997}). So one can attack the following question by first examining the obstruction of Lemma~\ref{lem0}, using algebraic methods. 
\\

\noindent \textbf{Question}: Is there a classification of L-space satellite operations?
\\

Another interesting direction one can pursue, encouraged by the computations done in this paper, is to calculate the Alexander polynomials $\Delta_{K}(T)$ of twisted $(p, q)$ torus knots with $q = kp \pm 1$ or more generally with $q$ an arbitrary non-zero integer. In \cite{Morton2006}, Morton gives a closed formula for $\Delta_{K}(T)$ where $K = K(p, q; 2, r)$ and $p>q>0$. 

\bibliographystyle{amsplain}

\bibliography{Reference}

\end{document}